\documentclass[11pt]{amsart} 

\usepackage{amssymb,tensor,enumitem,verbatim}
\usepackage{hyperref}
\usepackage[margin=1in]{geometry}
\usepackage[all,cmtip]{xy}
\usepackage[sort]{cite}
\usepackage{xcolor}

\newcommand{\mylabel}[2]{#2\def\@currentlabel{#2}\label{#1}}

\newcommand{\subgrp}[1]{\langle #1 \rangle}

\newcommand{\set}[1]{\{ #1 \}}

\newcommand{\abs}[1]{| #1 |}

\newcommand{\wh}[1]{\widehat{ #1}}
\newcommand{\wt}[1]{\widetilde{ #1}}
\newcommand{\ol}[1]{\overline{#1}}

\newcommand{\da}[1]{\!\!\downarrow_{#1}}

\DeclareMathOperator{\ad}{ad}

\DeclareMathOperator{\chr}{char}

\DeclareMathOperator{\Der}{Der}

\DeclareMathOperator{\Ext}{Ext}
\DeclareMathOperator{\flatdim}{flatdim}
\DeclareMathOperator{\Forget}{Forget}

\DeclareMathOperator{\gr}{gr}
\DeclareMathOperator{\opH}{H}
\DeclareMathOperator{\Hom}{Hom}
\DeclareMathOperator{\id}{id}

\DeclareMathOperator{\Lie}{Lie}
\DeclareMathOperator{\Proj}{Proj}
\DeclareMathOperator{\projdim}{projdim}
\DeclareMathOperator{\Max}{MaxSpec}
\DeclareMathOperator{\res}{res}
\DeclareMathOperator{\Spc}{Spc}
\DeclareMathOperator{\Spec}{Spec}
\DeclareMathOperator{\supp}{supp}

\newcommand{\stmodVL}{\stmod_{V(L)}}
\newcommand{\StModVL}{\StMod_{V(L)}}
\newcommand{\stmodVLbar}{\stmod_{V(\Lbar)}}
\newcommand{\StModVLbar}{\StMod_{V(\Lbar)}}
\newcommand{\stsmodVL}{\stsmod_{V(L)}}
\newcommand{\StsModVL}{\StsMod_{V(L)}}

\newcommand{\stmod}{\mathsf{stmod}}
\newcommand{\StMod}{\mathsf{StMod}}
\newcommand{\stsmod}{\mathsf{st\text{-}smod}}
\newcommand{\StsMod}{\mathsf{St\text{-}sMod}}

\newcommand{\Cbul}{C^\bullet}
\newcommand{\Hbul}{\opH^\bullet}

\newcommand{\Lbar}{\ol{L}}
\newcommand{\Ybar}{\ol{Y}}

\newcommand{\zero}{\ol{0}}
\newcommand{\one}{\ol{1}}

\newcommand{\ev}{\textup{ev}}
\newcommand{\odd}{\textup{odd}}

\newcommand{\ve}{\varepsilon}

\newcommand{\calJ}{\mathcal{J}}
\newcommand{\calN}{\mathcal{N}}
\newcommand{\calO}{\mathcal{O}}
\newcommand{\calP}{\mathcal{P}}
\newcommand{\calV}{\mathcal{V}}
\newcommand{\Chi}{\mathcal{X}}

\newcommand{\N}{\mathbb{N}}
\newcommand{\Z}{\mathbb{Z}}

\newcommand{\g}{\mathfrak{g}}

\newcommand{\fp}{\mathfrak{p}}
\newcommand{\fq}{\mathfrak{q}}

\newcommand{\wtg}{\wt{\g}}

\newcommand{\Lzero}{L_{\zero}}
\newcommand{\Lone}{L_{\one}}

\newcommand{\Vone}{V_{\ol{1}}}

\newcommand{\Vzero}{V_{\ol{0}}}

\newcommand{\sfmod}{\mathsf{mod}}
\newcommand{\sfMod}{\mathsf{Mod}}
\newcommand{\sfsmod}{\mathsf{smod}}
\newcommand{\sfsMod}{\mathsf{sMod}}
\newcommand{\sfsvec}{\mathsf{svec}}

\numberwithin{equation}{section}

\newtheorem{theorem}{Theorem}[section]
\newtheorem{proposition}[theorem]{Proposition}

\newtheorem{corollary}[theorem]{Corollary}
\newtheorem{lemma}[theorem]{Lemma}

\theoremstyle{definition}
\newtheorem{definition}[theorem]{Definition}

\newtheorem{example}[theorem]{Example}
\newtheorem{remark}[theorem]{Remark}

\title{Support varieties for Lie superalgebras in characteristic 2}

\author{Christopher M.\ Drupieski}
\address{Department of Mathematical Sciences,
		DePaul University,
		Chicago, IL 60614, USA}
\email{c.drupieski@depaul.edu}

\author{Jonathan R.\  Kujawa}
\address{Department of Mathematics \\
		University of Oklahoma \\
		Norman, OK 73019, USA}
\email{kujawa@math.ou.edu}

\thanks{The first author was supported in part by Simons Collaboration Grant for Mathematicians No.\ 4269055, and by a DePaul University Research Council paid faculty leave in Spring 2022. The second author was supported in part by Simons Collaboration Grant for Mathematicians No.\ 525043.}

\date{\today}


\allowdisplaybreaks

\begin{document}

\begin{abstract}
This paper investigates cohomology and support varieties for Lie superalgebras and restricted Lie superalgebras over a field of characteristic $2$.  The existence of an underlying ordinary Lie algebra allows us to obtain results that are still open in odd characteristic, and also to establish results that have no non-super analogues in characteristic $2$.
\end{abstract}

\dedicatory{Dedicated to Brian Parshall with gratitude.}

\maketitle


\section{Introduction}

\subsection{Overview}

In the mid-1980s, Friedlander and Parshall \cite{Friedlander:1987,Friedlander:1986b} introduced support varieties to the representation theory of restricted Lie algebras. Through their work and the work of those who followed, support varieties for restricted Lie algebras are comparatively well-understood. Almost four decades later, support varieties and their generalizations remain a subject of active investigation in non-semisimple representation theory, and the work of Friedlander and Parshall continues to serve as a model for those hoping to bring their methods to other settings. Relevant to this paper, one such setting is that of Lie superalgebras.

Roughly, a `super' version of a mathematical structure incorporates $\Z_{2}$-gradings on objects, together with $\pm$ signs that are inserted into the defining equations whenever elements of odd superdegree move past one another. Lie superalgebras arise in algebraic topology, spin representations of the symmetric group, mathematical physics, and elsewhere.  Work by Boe, Kujawa, and Nakano \cite{Boe:2010,Boe:2017}, the authors \cite{Drupieski:2019a,Drupieski:2021a,Drupieski:2021b}, Benson, Iyengar, Krause, and Pevtsova \cite{Benson:2021a,Benson:2022}, and others demonstrates that support varieties generalize in interesting and nontrivial ways to the super setting.

The $\pm$ signs that show up in super formulas are one of the defining features of the theory. However, when the ground field is of characteristic $2$ the $\pm$ signs effectively disappear, and by forgetting the $\Z_2$-gradings one seems to obtain just a classical, non-super version of the object in question. Indeed, it is sometimes assumed that super-mathematics over fields of characteristic $2$ is uninteresting.  Nevertheless, there is an extensive literature on the study of Lie superalgebras over fields of characteristic $2$; see, for example, \cite{Bouarroudj:2020} and the references therein.  The current paper hopefully adds additional evidence against that assumption.

The goal of this paper is to study support varieties for Lie superalgebras (both restricted and un-restricted) over fields of characteristic $2$. As we will see, the overlap between the super and non-super theories in characteristic $2$ is a feature, not a bug. One can use this overlap as leverage to prove results for the super-theory that are merely conjectural in odd characteristics, as well as results that have no non-super analogues in characteristic $2$. This paper thus exists at the intersection of topics that are well-understood (support varieties for restricted Lie algebras) and topics that are not (support varieties for graded objects).

\subsection{Main Results}

We begin in Section~\ref{S:Definitions} by recalling the definition of a Lie super\-algebra over a field of arbitrary characteristic. Roughly, a Lie superalgebra in characteristic $2$ is a $\Z_2$-graded ordinary Lie algebra $L = \Lzero \oplus \Lone$, together with a distinguished quadratic operator $q: \Lone \to \Lzero$. The quadratic operator effectively serves as a replacement for the map given in other characteristics by $x \mapsto \frac{1}{2}[x,x]$; see Remark \ref{rem:oddbracket}. A representation of a Lie superalgebra in characteristic $2$ is then a $\Z_2$-graded representation of the underlying ordinary Lie algebra $\Lbar$ that is `partially restricted' to ensure compatibility with the quadratic operator.

A restricted Lie superalgebra is a Lie superalgebra $L$ equipped with a $p$-operation $x \mapsto x^{[p]}$ on $\Lzero$ that makes $\Lzero$ into an ordinary restricted Lie algebra, and that makes $\Lone$ into a restricted $\Lzero$-module. In Lemma \ref{lemma:RLSA=RLA} we observe that, given a restricted Lie superalgebra over a field of characteristic $p=2$, the $2$-operation $(-)^{[2]}: \Lzero \to \Lzero$ and the quadratic operator $q: \Lone \to \Lzero$ can be combined to define a $2$-operation $z \mapsto z^{\set{2}}$ on all of $L$,
	\[
	z = z_{\zero} + z_{\one} \quad \mapsto \quad z^{\set{2}} := (z_{\zero})^{[2]} + q(z_{\one}) + [z_{\one},z_{\zero}],
	\]
which then makes the underlying ordinary Lie algebra $\Lbar$ of $L$ into a restricted Lie algebra.\footnote{This observation has doubtlessly been made before and is presumably well-known to characteristic $2$ experts.} Moreover, this restricted structure on $\Lbar$ then induces a canonical identification between $V(L)$, the restricted enveloping superalgebra of the restricted Lie superalgebra $L$, and $V(\Lbar)$, the restricted enveloping algebra of the ordinary restricted Lie algebra $\Lbar$.

\subsubsection{Restricted Lie superalgebras}

In Section \ref{S:support-for-restricted} we explore in detail the implications for support varieties of the canonical identification $V(L) = V(\Lbar)$, for $L$ a finite-dimensional restricted Lie superalgebra over a (for simplicity, algebraically closed) field $k$ of characteristic $2$. The cohomology ring $\Hbul (V(L), k)$ is a finitely-generated commutative $k$-algebra that is naturally graded by $\Z \times \Z_{2}$; the first component of the grading is the `external' cohomological grading and the second component is the `internal' grading induced by the superalgebra structure of $V(L)$. Write $\abs{V(L)}$ for the maximal ideal spectrum of the $k$-algebra $\Hbul(V(L),k)$. We show in \eqref{eq:V(L)spectrum} that $\abs{V(L)}$ identifies with the restricted nullcone of $\Lbar$:
	\[
	\abs{V(L)} \simeq \calN_1(\Lbar) = \set{ z \in \Lbar: z^{\set{2}} = 0} = \set{ z_{\zero}+z_{\one} \in L : (z_{\zero})^{[2]} + q(z_{\one}) = 0 \text{ and } [z_{\one},z_{\zero}] = 0}.
	\]
This extends to characteristic $p=2$ a description previously obtained by the authors for restricted Lie superalgebras in characteristics $p \geq 3$; cf.\ \cite[Proposition 5.4.4]{Drupieski:2019a}.

Given a finite-dimensional $V(L)$-supermodule $M$, let $I(M)$ be the annihilator ideal for the cup product action of $\Hbul(V(L),k)$ on $\Ext_{V(L)}^\bullet(M,M)$. The support variety of $M$, denoted $\abs{V(L)}_M$, is the Zariski closed subset of $\abs{V(L)}$ defined by the ideal $I(M)$. Each nonzero element $z \in \calN_1(\Lbar)$ generates a $k$-subalgebra of $V(L) = V(\Lbar)$ of the form $k[z]/(z^2)$. Appealing to support variety theory for ordinary restricted Lie algebras, the support variety $\abs{V(L)}_M$ then admits the following ``rank variety'' description:
	\begin{equation}\label{E:rankvariety}
	\abs{V(L)}_{M} \simeq \set{ z \in \calN_1(\Lbar) : M \text{ is not free as a $k[z]/(z^2)$-module}} \cup \set{0}.
	\end{equation}
This immediately implies that support varieties for $V(L)$-supermodules satisfy the Tensor Product Property (Proposition~\ref{prop:tensor-product-property}). In Proposition~\ref{prop:rank-projdim}, we reformulate the condition that $M$ is not free as a $k[z]/(z^2)$-module in terms of the projective dimension of a pullback of $M$ to a hyper\-surface ring. As discussed in Remark \ref{rem:detecting-finite=proj-dim}, this alternate characterization of $\abs{V(L)}_M$ can be interpreted as confirming, for characteristic $2$, a conjecture we previously made for the support varieties of restricted Lie superalgebras (and more generally, for arbitrary infinitesimal supergroup schemes) over fields of odd characteristic.

Now let $\stsmodVL$ be the stable category of finite-dimensional $V(L)$-supermodules.  The Hopf algebra structure on $V(L)$ makes $\stsmodVL$ into a tensor triangulated category. Following Balmer, one can use the tensor triangulated structure to define the spectrum $\Spc (\stsmodVL )$ of $\stsmodVL$ and the support function 
	\[
	\supp' : \stsmodVL \to 2^{\Spc (\stsmodVL )}.
	\]
Here $2^{\Spc (\stsmodVL )}$ denotes the power set of $\Spc (\stsmodVL )$.  Given its abstract definition, it is not surprising that it is desirable to have a more concrete description of the Balmer spectrum and support function.

Let $X_s = \Proj_{s}(A)$ denote the set of all bi-homogeneous ideals in $A = \Hbul (V(L),k)$ (i.e., ideals that are homogeneous with respect to the $\Z \times \Z_{2}$-grading) that are prime among the bi-homogeneous ideals, and which are properly contained in the maximal ideal $\opH^{> 0}(V(L),k)$ of elements of positive cohomological degree. We view $X_s$ as a topological space via the Zariski topology. We note that $X_s$ is a coarser space than the usual projective spectrum $\Proj(A)$ (which consists of prime ideals in $\Hbul(V(L),k)$ which are homogeneous with respect to just the $\Z$-component of the grading).

In Section~\ref{subsec:tensor-triangular-geometry} we show that $X_s$ is a Zariski space and that there is a support function 
	\[
	\supp_{s} : \stsmodVL \to 2^{X_s}.
	\]
Applying work of Dell'Ambrogio \cite{DellAmbrogio:2010} and of Benson, Iyengar, Krause, and Pevtsova \cite{Benson:2018}, we show in Theorem \ref{T:BalmerSpectrum} that there is a canonical homeomorphism
	\[
	\Proj_{s}(A) \simeq \Spc (\stsmodVL )
	\]
that identifies the two support functions $\supp'$ and $\supp_s$. This gives a non-categorical realization of the Balmer spectrum and support function. We further deduce that the pair $(\Proj_{s}(A), \supp_{s} )$ provides a classification of the thick tensor ideal subcategories of $\stsmodVL$. Concretely, this implies that if $M$ and $N$ are two finite-dimensional $V(L)$-supermodules, then $M$ can be obtained from $N$ by a finite sequence of the standard operations available in $\stsmodVL$ (direct sums and summands, tensor products, extensions, syzygies) if and only if $\abs{V(L)}_{M} \subseteq \abs{V(L)}_{N}$. Thus the question of whether one supermodule can be built from another can be reformulated into a statement about their support varieties, which in turn is reduced to a linear-algebraic problem by \eqref{E:rankvariety}.

\subsubsection{Lie superalgebras}

In Section~\ref{S:support-for-non-restricted} we turn to the case of an arbitrary finite-dimensional Lie superalgebra $L$ over $k$. Again, the cohomology ring $\Hbul (L,k) = \Hbul(U(L),k)$ is a finitely-generated commutative $k$-algebra graded by $\Z \times \Z_{2}$, and one can define support varieties for finite-dimensional $L$-supermodules. Writing $\abs{U(L)}$ for the maximal ideal spectrum of $\Hbul(L,k)$, we show in Theorem \ref{theorem:spectrum-odd-nullcone} that $\abs{U(L)}$ identifies with the odd nullcone of $L$:
	\[
	\abs{U(L)} \simeq \calN_{\odd}(L) := \set{ x \in \Lone : q(x) = 0}.\footnote{There is also a Frobenius twist involved in the identification, which we ignore for the sake of the introduction.}
	\]
Likewise, Theorem~\ref{theorem:support=rank} shows that the support variety of a finite-dimensional $L$-supermodule $M$ admits a rank variety description:
	\[
	\abs{U(L)}_M \simeq \set{x \in \calN_{\odd}(L) : M \text{ is not free as a $k[x]/(x^2)$-supermodule}} \cup \set{0}.
	\]
This description implies that the support varieties of finite-dimensional $L$-supermodules satisfy the Tensor Product Property (Corollary \ref{cor:tensor-product-property-U(L)}). Then Theorem~\ref{T:DetectFiniteProjDimension} shows that $\abs{U(L)}_M = \set{0}$ if and only if $M$ has finite projective dimension as a $U(L)$-supermodule. In particular, $\calN_{\odd}(L) = \set{0}$ if and only if $U(L)$ is of finite global dimension. This extends to $p=2$ a result due of B{\o}gvad \cite{Bo-gvad:1984} in characteristic zero (see also \cite[Chapter 17]{Musson:2012}), and due to the authors in characteristic $p \geq 3$ \cite{Drupieski:2021b}.

\subsection{Future directions}

Since Premet's proof of the Kac--Weisfieler conjecture \cite{Premet:1995}, $p$-divisibility of a module's dimension has been known to be closely related to the dimension of the module's support variety. Some divisibility results along these lines are known in super-representation theory (see, e.g., \cite{Boe:2009,Wang:2009}), but much remains unknown, however, and it would be interesting to improve upon these results.  

It would also be interesting to compute, in all positive characteristics, the Balmer spectrum for the stable module category of finite-dimensional $V(L)$-supermodules, and for the singularity category of finite-dimensional $U(L)$-supermodules. More generally, there are conjectures and questions for Lie superalgebras and infinitesimal supergroup schemes that the present paper addresses, but which remain open in odd characteristics.  For $V(L)$ we were able to take advantage of coincidences that happen only in characteristic $2$. In odd characteristic new ideas will be required for both superalgebras.  

Finally, interesting alternate notions of super-mathematics in characteristic $2$ were proposed in \cite{Venkatesh:2016,Kaufer:2018}. It would be interesting to consider the questions covered in this paper within that framework.

\subsection{Conventions}

Throughout, $k$ will be a field of positive characteristic $p$ (usually $p=2$). All vector spaces will be $k$-vector spaces, and all unadorned tensor products will denote tensor products over $k$. Let $\N = \set{0,1,2,3,\ldots}$ be the set of non-negative integers. 

Set $\Z_2 = \Z/2\Z = \{ \zero,\one \}$. Following the literature, we use the prefix `super' to indicate that an object is $\Z_2$-graded. We denote the decomposition of a vector superspace into its $\Z_2$-homogeneous components by $V = \Vzero \oplus \Vone$, calling $\Vzero$ and $\Vone$ the even and odd subspaces of $V$, respectively, and writing $\ol{v} \in \Z_2$ to denote the superdegree of a homogeneous element $v \in \Vzero \cup \Vone$. When written without additional adornment, we consider the field $k$ to be a superspace concentrated in even super\-degree. Whenever we state a formula in which homogeneous degrees of elements are specified, we mean that the formula is true as written for homogeneous elements, and that it extends linearly to non-homogeneous elements.

\subsection{Acknowledgments}

The authors thank Luchezar Avramov for prodding us to address what happens in characteristic two.

\section{Definitions}\label{S:Definitions}

\subsection{Lie superalgebras} \label{subsec:lie-superalgebras}

\begin{definition} \label{def:Liesuperalgebra}
Let $k$ be a field (of any characteristic). A Lie superalgebra over $k$ is a $\Z_2$-graded $k$-vector space $L = \Lzero \oplus \Lone$, together with an even (i.e., $\Z_2$-degree preserving) bilinear map $[ \, , \,]: L \otimes L \to L$ and a quadratic operator $q: \Lone \to \Lzero$, such that the following axioms hold:
	\begin{enumerate}[label={(L\arabic*)}]
	\item \label{item:alternating} $[x,y] = -(-1)^{\ol{x} \cdot \ol{y}} [y,x]$ for all homogeneous $x,y \in L$.

	\item \label{item:Jacobiidentity} $[x,[y,z]] = [[x,y],z] + (-1)^{\ol{x} \cdot \ol{y}} [y,[x,z]]$ for all homogeneous $x,y,z \in L$ (Jacobi identity).

	\item \label{item:evenselfcommute} $[x,x] = 0$ for all $x \in \Lzero$.

	\item \label{item:oddcommutewithsquare} $[y,[y,y]] = 0$ for all $y \in \Lone$.

	\item \label{item:oddbracketequalsq} $[x,y] = q(x+y) - q(x) - q(y)$ for all $x,y \in \Lone$.

	\item \label{item:oddadjointtoq} $[y,[y,z]] = [q(y),z]$ for all $y \in \Lone$ and $z \in L$.
	\end{enumerate}
\end{definition}

\begin{remark} \label{rem:oddbracket}
Let $L$ be a Lie superalgebra, and let $y \in \Lone$. Then \ref{item:oddbracketequalsq} asserts that
	\begin{equation} \label{eq:oddbracket}
	[y,y] = q(y+y) - q(y) - q(y) = q(2 \cdot y) - 2 \cdot q(y) = 4 \cdot q(y) - 2  \cdot q(y) = 2 \cdot q(y).
	\end{equation}
For $p \neq 2$, this shows that $q(y) = \frac{1}{2}[y,y]$, while for $p = 2$, it shows that $[y,y] = 0$ for all $y \in \Lone$. Thus for $p=2$, if we forget the quadratic operator $q: \Lone \to \Lzero$, the Lie bracket $[\, , \,]: L \otimes L \to L$ makes $L$ into an ordinary Lie algebra, which we may denote $\Lbar$.
\end{remark}

\begin{remark} \label{rem:redundancies}
Depending on the characteristic $p$ of the field $k$, some of the axioms of a Lie superalgebra are redundant:
\begin{enumerate}
\item If $p \neq 2$, then \ref{item:evenselfcommute} is a consequence of \ref{item:alternating}.

\item If $p \neq 3$, then \ref{item:oddcommutewithsquare} is a consequence of \ref{item:alternating} and \ref{item:Jacobiidentity}.

\item Suppose $p \neq 2$. Then \ref{item:oddbracketequalsq} and \ref{item:oddadjointtoq} follow from \ref{item:alternating}, \ref{item:Jacobiidentity}, and the formula $q(y) = \frac{1}{2}[y,y]$. Conversely, if $L$ is a $k$-superspace equipped with an even bilinear map $[\, , \,]: L \otimes L \to L$ satisfying \ref{item:alternating}--\ref{item:oddcommutewithsquare}, then $q: \Lone \to \Lzero$ defined by $q(y) = \frac{1}{2}[y,y]$ satisfies \ref{item:oddbracketequalsq} and \ref{item:oddadjointtoq}.
\end{enumerate}
It is necessary to assume \ref{item:evenselfcommute} when $p = 2$, and to assume \ref{item:oddcommutewithsquare} when $p=3$, to ensure that $L$ can be embedded into its universal enveloping algebra $U(L)$; see \cite[1.1.10 and 3.2.3]{Bahturin:1992}.
\end{remark}

\begin{remark}
For those who (like the authors of this paper) have spent most of their lives only considering Lie superalgebras over fields of characteristic $p \neq 2$, the inclusion of the quadratic operator $q: \Lone \to \Lzero$ in Definition \ref{def:Liesuperalgebra} may require explanation. Remarks \ref{rem:oddbracket} and \ref{rem:redundancies} provide some motivation for why it is natural to include the quadratic operator in the definition of a Lie super\-algebra. Another source of motivation is the fact that it permits a uniform, characteristic-free description of the Koszul resolution (Theorem \ref{theorem:Koszulresolution}). Here is yet another motivating observation:

Let $k$ be a field (of any characteristic), let $A$ be a $k$-superalgebra, and let $D \in \Der_k(A,A)_{\one}$ be an odd superderivation of $A$. Thus $D$ is an odd linear map such that for all homogeneous elements $a,b \in A$, one has $D(ab) = D(a) \cdot b + (-1)^{\ol{a}} a \cdot D(b)$.
Then
	\begin{align*}
	D^2(ab) &= D\left( D(a) \cdot b + (-1)^{\ol{a}} a \cdot D(b) \right) \\
	&= D^2(a) \cdot b + (-1)^{\ol{a} + \one} D(a) \cdot D(b) + (-1)^{\ol{a}} D(a) \cdot D(b) + a \cdot D^2(b) \\
	&= D^2(a) \cdot b + a \cdot D^2(b),
	\end{align*}
so $D^2 \in \Der_k(A,A)_{\zero}$. Now, one of the ways to describe the Lie superalgebra $\Lie(G)$ of an affine algebraic supergroup $G$ is as the space of left-invariant superderivations of its coordinate super\-algebra $k[G]$.\footnote{Recall that, writing $\Delta: k[G] \to k[G] \otimes k[G]$ for the coproduct on $k[G]$, a superderivation $D: k[G] \to k[G]$ is left-invariant if $\Delta \circ D = (\id \otimes D) \circ \Delta$.} The preceding calculation then implies that, in arbitrary characteristic, the squaring map $D \mapsto D^2$ defines a quadratic operator $q: \Lie(G)_{\one} \to \Lie(G)_{\zero}$.
\end{remark}

\begin{example}
Let $k$ be a field (of any characteristic), and let $A$ be a $k$-superalgebra. Define an even bilinear map $[\, , \,]: A \otimes A \to A$ and a quadratic operator $q: A_{\one} \to A_{\zero}$ by
	\begin{itemize}
	\item $[a,b] = ab - (-1)^{\ol{a} \cdot \ol{b}}  ba$ for $a,b \in A$ homogeneous, and
	\item $q(a) = a^2$ for $a \in A_{\one}$.
	\end{itemize}
These operations define a Lie superalgebra structure on the underlying $k$-superspace of $A$, which we may denote by $\Lie(A)$.
\end{example}

\begin{definition} \label{def:envelopingalgebra}
Let $k$ be a field (of any characteristic), let $L$ be a Lie superalgebra over $k$, and let $T(L) = \bigoplus_{i \geq 0} L^{\otimes i}$ be the tensor (super)algebra on $L$. The universal enveloping algebra of $L$ is the $k$-superalgebra $U(L) := T(L)/I$, where $I$ is the ideal in $T(L)$ generated by the following elements:
	\begin{itemize}
	\item $x \otimes y - (-1)^{\ol{x} \cdot \ol{y}} y \otimes x - [x,y]$ for $x,y \in L$ homogeneous, and
	\item $y \otimes y - q(y)$ for $y \in \Lone$.
	\end{itemize}
\end{definition}

The proof of the following theorem is standard.

\begin{theorem}[Poincar\'{e}--Birkhoff--Witt Theorem for $U(L)$]
Let $L$ be a Lie superalgebra over a field $k$ (of any characteristic), and let $X = \set{x_\lambda : \lambda \in \Omega}$ be an ordered basis for $L$. Then $U(L)$ admits a $k$-basis consisting of all monomials of the form $x_1 x_2 \cdots x_n$, where $n \in \N$, $x_i \in X$, $x_i \leq x_{i+1}$, and $x_i \neq x_{i+1}$ if $x_i \in \Lone$. In particular, the inclusion $L = L^{\otimes 1} \hookrightarrow T(L)$ induces an injective Lie super\-algebra homomorphism $\iota: L \hookrightarrow \Lie(U(L))$.
\end{theorem}

\begin{remark}
Let $L$ be a Lie superalgebra over a field $k$ of characteristic $p$.
	\begin{enumerate}
	\item The algebra $U(L)$ is a cocommutative Hopf superalgebra, with the subspace $L$ consisting of primitive elements. If $p = 2$, then $U(L)$ is a Hopf algebra in the ordinary sense.
	\item Suppose $p = 2$, and let $\Lbar$ be the ordinary ($\Z_2$-graded) Lie algebra obtained by forgetting the quadratic operator $q: \Lone \to \Lzero$. Then
		\[
		U(L) = U(\Lbar)/\subgrp{y^2 - q(y) : y \in \Lone},
		\]
	and the canonical quotient map $U(\Lbar) \twoheadrightarrow U(L)$ is a Hopf (super)algebra homomorphism.
	\item Suppose $p=2$. More-or-less by definition, a (left) $L$-supermodule is a (left) supermodule over the universal enveloping superalgebra $U(L)$. Because $y^2 - q(y) = 0$ in $U(L)$ for each $y \in \Lone$, but $y^2 - q(y) \neq 0$ in $U(\Lbar)$, this shows that the category of $L$-supermodules is not simply the category of ($\Z_2$-graded) $\Lbar$-modules.
	\end{enumerate}
\end{remark}

\subsection{Restricted Lie superalgebras} \label{subsec:restricted-lie-superalgebras}

\begin{definition}
Let $k$ be a field of characteristic $p \geq 2$, and let $L$ be a Lie superalgebra over $k$. We say that $L$ is a restricted Lie superalgebra if there exists a map $x \mapsto x^{[p]}$ on $\Lzero$ (called a \emph{$p$-map} or a \emph{$p$-operation}) such that the following are satisfied for all $x,y \in \Lzero$:
	\begin{enumerate}[label={(R\arabic*)}]
	\item \label{item:p-homogeneous} $(\alpha \cdot x)^{[p]} = \alpha^p \cdot x^{[p]}$ for all $\alpha \in k$,
	\item \label{item:ad-pth-power} $\ad(x^{[p]}) = \ad(x)^p$
	\item \label{item:p-on-sums} $(x+y)^{[p]} = x^{[p]} + y^{[p]} + \sum_{i=1}^{p-1} s_i(x,y)$, where $i \cdot s_i(x,y)$ is the coefficient of $t^{i-1}$ in the formal expression $\ad(t \cdot x + y)^{p-1}(x)$.
	\end{enumerate}
In other words, $L$ is a restricted Lie superalgebra if $\Lzero$ is a restricted Lie algebra in the usual sense, and $\Lone$ is a restricted $\Lzero$-module.
\end{definition}

\begin{remark}
For $p=2$, the condition \ref{item:p-on-sums} becomes $(x+y)^{[2]} = x^{[2]} + y^{[2]} + [y,x]$.
\end{remark}

\begin{example}
Let $k$ be a field of characteristic $p \geq 2$, and let $A$ be a $k$-superalgebra. Then the map $a \mapsto a^p$ defines a $p$-operation on $A_{\zero}$ that makes $\Lie(A)$ into a restricted Lie superalgebra.
\end{example}

\begin{definition}
Let $k$ be a field of characteristic $p \geq 2$, and let $L$ be a restricted Lie superalgebra over $k$. The restricted enveloping algebra of $L$ is the $k$-superalgebra
	\[
	V(L) := U(L) / \subgrp{ x^p - x^{[p]} : x \in \Lzero}.
	\]
\end{definition}

The proof of the following theorem is also standard.

\begin{theorem}[Poincar\'{e}--Birkhoff--Witt Theorem for $V(L)$] \label{thm:PBW-restricted}
Let $L$ be a restricted Lie superalgebra over a field $k$ of characteristic $p \geq 2$, and let $X = \set{x_\lambda : \lambda \in \Omega}$ be an ordered basis for $L$. Then $V(L)$ admits a $k$-basis consisting of all monomials of the form $x_1^{\lambda_1} x_2^{\lambda_2}  \cdots x_n^{\lambda_n}$, where $n \in \N$, $x_i \in X$, $x_1 < x_2 < \cdots < x_n$, $0 \leq \lambda_i < p$ if $x_i \in \Lzero$, and $0 \leq \lambda_i < 2$ if $x_i \in \Lone$. The natural inclusion $L \hookrightarrow T(L)$ induces an injective homomorphism of restricted Lie superalgebras $\iota: L \hookrightarrow \Lie(V(L))$.
\end{theorem}

Let $G$ be an abelian group, written additively. Recall that a \emph{$G$-graded Lie algebra} is a (ordinary) Lie algebra $L$ together with a family of subspaces $(L_g)_{g \in G}$ of $L$ such that $L = \bigoplus_{g \in G} L_g$ and $[L_g,L_h] \subseteq L_{g+h}$ for all $g,h \in G$. If $L$ is a (ordinary) restricted Lie algebra with $p$-map $x \mapsto x^{[p]}$, then $L$ is a \emph{restricted $G$-graded Lie algebra} if the $p$-map satisfies $(L_g)^{[p]} \subseteq L_{pg}$ for all $g \in G$.

\begin{lemma} \label{lemma:RLSA=RLA}
Let $L$ be a restricted Lie superalgebra over a field $k$ of characteristic $p = 2$, with $p$-operation $x \mapsto x^{[2]}$ on $\Lzero$, and let $\Lbar$ be the ordinary Lie algebra over $k$ obtained by forgetting the quadratic operator $q: \Lone \to \Lzero$. Given $z \in L$, let $z = z_{\zero} + z_{\one}$ be the decomposition of $z$ into its even and odd parts. Then the map
	\begin{equation} \label{eq:2-operation}
	z \mapsto z^{\set{2}} := (z_{\zero})^{[2]} + q(z_{\one}) + [z_{\one},z_{\zero}]
	\end{equation}
defines a $p$-operation on $\Lbar$ that makes $\Lbar$ into a restricted ($\Z_2$-graded) Lie algebra over $k$. The canonical quotient map $U(\Lbar) \twoheadrightarrow U(L)$ induces a Hopf (super)algebra isomorphism $V(\Lbar) \cong V(L)$.
\end{lemma}

\begin{proof}
Let $x, y \in \Lbar$, and let $x = x_{\zero} + x_{\one}$ and $y = y_{\zero} + y_{\one}$ be the decompositions of $x$ and $y$ into their $\Z_2$-homogeneous components. First, for $\alpha \in k$, one has
	\begin{align*}
	(\alpha \cdot x)^{\set{2}} &= (\alpha \cdot x_{\zero} + \alpha \cdot x_{\one})^{\set{2}} \\
	&= (\alpha \cdot x_{\zero})^{[2]} + q(\alpha \cdot x_{\one}) + [\alpha \cdot x_{\one},\alpha \cdot x_{\zero}] \\
	&= \alpha^2 \cdot (x_{\zero})^{[2]} + \alpha^2 \cdot q(x_{\one}) + \alpha^2 \cdot [x_{\one},x_{\zero}] \\
	&= \alpha^2 \cdot x^{\set{2}},
	\end{align*}
so \ref{item:p-homogeneous} is satisfied by the map $z \mapsto z^{\set{2}}$. Second,
	\begin{align*}
	\ad(x^{\set{2}})(y) &= [x^{\set{2}},y] \\
	&= [(x_{\zero})^{[2]} + q(x_{\one}) + [x_{\one},x_{\zero}], y] \\
	&= [(x_{\zero})^{[2]},y] + [q(x_{\one}),y] + [[x_{\one},x_{\zero}], y] \\
	&= [x_{\zero},[x_{\zero},y]] + [x_{\one},[x_{\one},y]] + [x_{\one},[x_{\zero}, y]] + [x_{\zero},[x_{\one},y]] \\
	&= [x_{\zero}+x_{\one},[x_{\zero}+x_{\one},y]] \\
	&= \ad(x_{\zero}+x_{\one})^2(y), 
	\end{align*}
so \ref{item:ad-pth-power} is satisfied by the map $z \mapsto z^{\set{2}}$. Third,
	\begin{align*}
	&(x+y)^{\set{2}} \\
	&= ( (x_{\zero}+y_{\zero}) + (x_{\one}+y_{\one}))^{\set{2}} \\
	&= (x_{\zero}+y_{\zero})^{[2]} + q(x_{\one}+y_{\one}) + [x_{\one}+y_{\one}, x_{\zero}+y_{\zero}] \\
	&= (x_{\zero})^{[2]} + (y_{\zero})^{[2]} + [y_{\zero},x_{\zero}] + [x_{\one},y_{\one}]+q(x_{\one}) + q(y_{\one}) + [x_{\one},x_{\zero}] + [x_{\one},y_{\zero}] + [y_{\one}, x_{\zero}] + [y_{\one},y_{\zero}] \\
	&= x^{\set{2}} + y^{\set{2}} + [y_{\zero},x_{\zero}] + [y_{\one},x_{\one}] + [y_{\zero},x_{\one}] + [y_{\one}, x_{\zero}] \\
	&= x^{\set{2}} + y^{\set{2}} + [y_{\zero}+y_{\one},x_{\zero}+x_{\one}] \\
	&= x^{\set{2}} + y^{\set{2}} + [y,x]
	\end{align*}
so \ref{item:p-on-sums} is satisfied by the map $z \mapsto z^{\set{2}}$. Finally, let $j: U(\Lbar) \twoheadrightarrow U(L)$ and $\pi: U(L) \twoheadrightarrow V(L)$ be the canonical quotient maps. For $z = z_{\zero} + z_{\one} \in L$, in the algebra $U(\Lbar)$ one has
	\begin{align*}
	z^2 - z^{\set{2}} &= (z_{\zero}+z_{\one})^2 - \big( (z_{\zero})^{[2]} + q(z_1) + [z_{\one},z_{\zero}] \big) \\
	&= \big( (z_{\zero})^2 + z_{\zero}z_{\one} + z_{\one}z_{\zero} + (z_{\one})^2 \big) - \big( (z_{\zero})^{[2]} + q(z_1) + [z_{\one},z_{\zero}] \big) \\
	&= \big( (z_{\zero})^2 - (z_{\zero})^{[2]} \big) + \big( (z_{\one})^2 - q(z_1) \big).
	\end{align*}
This is a sum of a defining relation in $U(L)$ and a defining relation in $V(L)$, so the composite algebra homomorphism $\pi \circ j: U(\Lbar) \twoheadrightarrow V(L)$ factors through a surjective homomorphism $V(\Lbar) \twoheadrightarrow V(L)$, which can then be seen to be an injection (and hence an isomorphism) by Theorem \ref{thm:PBW-restricted}.
\end{proof}

\begin{remark}
The converse of Lemma \ref{lemma:RLSA=RLA} is true. That is, let $\Lbar$ be a restricted $\Z_2$-graded Lie algebra over a field of characteristic $p=2$, with $2$-map $z \mapsto z^{\set{2}}$. For $z \in \Lbar_{\one}$, set $q(z) = z^{\set{2}} \in \Lbar_{2 \cdot \one} = \Lbar_{\zero}$, and for $z \in \Lbar_{\zero}$ set $z^{[2]} = z^{\set{2}} \in \Lbar_{2 \cdot \zero} = \Lbar_{\zero}$. Then $q: \Lbar_{\one} \to \Lbar_{\zero}$ defines a quadratic operator that makes the underlying $\Z_2$-graded Lie algebra of $\Lbar$ into a Lie superalgebra, which we may denote $L$, and the map $z \mapsto z^{[2]}$ defines a $p$-operation on $\Lzero$ that makes $L$ into a restricted Lie superalgebra.
\end{remark}

\section{Support for restricted Lie superalgebras in characteristic 2} \label{S:support-for-restricted}

\subsection{Cohomology} \label{subsec:restricted-cohomology}

Let $L$ be a finite-dimensional restricted Lie super\-algebra over a field $k$ of characteristic $p=2$, and let $V(L)$ be its restricted enveloping algebra. Let $\sfsMod = \sfsMod_{V(L)}$ be the category of left $V(L)$-super\-modules, and let $\sfsmod = \sfsmod_{V(L)}$ be the subcategory of finite-dimensional $V(L)$-supermodules. An object in $\sfsMod_{V(L)}$ is a $k$-superspace $M$ equipped with an even bilinear map $V(L) \otimes M \to M$ that makes $M$ into a left $V(L)$-module. Given $V(L)$-supermodules $M$ and $N$, one has
	\[
	\Hom_{\sfsMod}(M,N) = \Hom_{\sfsMod}(M,N)_{\zero} \oplus \Hom_{\sfsMod}(M,N)_{\one}
	\]
where $\Hom_{\sfsMod}(M,N)_{\zero}$ (resp.\ $\Hom_{\sfsMod}(M,N)_{\one}$) is the subspace of all even (resp.\ odd) $k$-linear maps $f: M \to N$ such that $f(u.m) = u.f(m)$ for all $u \in V(L)$ and $m \in M$.\footnote{If $p \neq 2$, then $\Hom_{\sfsMod}(M,N)_{\one}$ consists of all odd $k$-linear maps $f: M \to N$ such that $f(u.m) = (-1)^{\ol{u}} u.f(m)$ for all homogeneous $u \in V(L)$ and all $m \in M$.} The category $\sfsMod_{V(L)}$ is not an abelian category (because the kernel or image of a non-homogeneous homomorphism need not be the direct sum of its even and odd subspaces), but the underlying even subcategory $(\sfsMod_{V(L)})_{\ev}$ of $\sfsMod_{V(L)}$, having the same objects as $\sfsMod_{V(L)}$ but only the even homomorphisms as morphisms, is an abelian category.

\begin{remark}
For $p \neq 2$, the category $(\sfsMod_{V(L)})_{\ev}$ identifies with the category of left modules for the smash product algebra $V(L) \# k\Z_2$. Here $V(L)$ is considered as a $k\Z_2$-module algebra, with the nontrivial element of $\Z_2$ acting on $V(L)$ via the sign automorphism $u \mapsto (-1)^{\ol{u}} u$. This identification of categories no longer holds for $p=2$, because then $k\Z_2$-modules no longer decompose into $\pm 1$ eigenspaces, and hence you can no longer use the action of $k\Z_2$ to recover the superspace structure of modules, and you cannot use $k\Z_2$-equivariance of maps to ensure that you are only dealing with even $V(L)$-supermodule homomorphisms.
\end{remark}

Each $V(L)$-supermodule is the image, via an even homomorphism, of a direct sum of copies of $V(L)$ and its parity shift $\Pi(V(L))$.\footnote{For $p=2$, the parity shift $\Pi(M)$ of a $V(L)$-supermodule $M$ is the $V(L)$-supermodule that has the same underlying $k$-vector space as $M$, and the same action of $V(L)$ on that $k$-vector space, but with $\Pi(M)_{\zero} = M_{\one}$ and $\Pi(M)_{\one} = M_{\zero}$. For $p \neq 2$, there is more than once choice of parity change functor; see \cite[\S2.3]{Drupieski:2021a}.} Then the category $(\sfsMod_{V(L)})_{\ev}$ has enough projectives, so given $V(L)$-supermodules $M$ and $N$, we can define $\Ext_{\sfsMod}^n(M,N)$ to be the $n$-th right derived functor of the contravariant functor $\Hom_{\sfsMod}(-,N) : (\sfsMod_{V(L)})_{\ev} \to \sfsvec_k$, i.e.,
	\[
	\Ext_{\sfsMod_{V(L)}}^n(M,N) = \big( R^n \Hom_{\sfsMod_{V(L)}}(-,N) \big)(M).
	\]
The category $(\sfsMod_{V(L)})_{\ev}$ also has enough injectives: Given a $V(L)$-supermodule $N$, define $\varphi: N \to \Hom_k(V(L),N)$ by $\varphi(n) = \varphi_n$, where $\varphi_n(u) = (-1)^{\ol{u} \cdot \ol{n}} u.n$. Then $\varphi$ is an even $V(L)$-super\-module homomorphism, where the left $V(L)$-module structure on $\Hom_k(V(L),N)$  comes from the right action of $V(L)$ on itself, and the Hom-tensor adjunction
	\[
	\Hom_{V(L)}(M,\Hom_k(V(L),N)) \cong \Hom_k(V(L) \otimes_{V(L)} M, N) \cong \Hom_k(M,N)
	\]
implies that $\Hom_k(V(L),N)$ is injective in $(\sfsMod_{V(L)})_{\ev}$. In fact, $V(L)$ is itself an injective object in $(\sfsMod_{V(L)})_{\ev}$. This can be seen from \eqref{eq:Homsmod=Hommod}, using the fact that $V(L) = V(\Lbar)$ is self-injective in the category of arbitrary $V(\Lbar)$-modules because $V(\Lbar)$ is a finite-dimensional Hopf algebra; cf.\ \cite[Lemma 2.3.2]{Drupieski:2019a} for $p \neq 2$. From this it follows that $(\sfsMod_{V(L)})_{\ev}$ is a Frobenius category.

Let $\sfMod = \sfMod_{V(L)}$ be the category of arbitrary (i.e., not necessarily $\Z_2$-graded) $V(L)$-modules. The forgetful functor, $\Forget: \sfsMod_{V(L)} \to \sfMod_{V(L)}$, which simply ignores the $\Z_2$-gradings on objects, defines a canonical identification
	\begin{equation} \label{eq:Homsmod=Hommod}
	\Hom_{\sfsMod}(M,N) = \Hom_{\sfMod}(M,N) \quad \text{for all $M,N \in \sfsMod_{V(L)}$.}\footnote{In general, the forgetful functor acts linearly on morphisms, with $\Forget(f) = f$ if $f$ is even, and $\Forget(f) = \wt{f}$ if $f$ is odd. Here $\wt{f}: M \to N$ is the $k$-linear map defined by $\wt{f}(m) = (-1)^{\ol{m}} \cdot f(m)$. For $p = 2$, this means that $\Forget(f) = f$ for all $f \in \Hom_{\sfsMod}(M,N)$.}
	\end{equation}
More generally, since each object in $\sfsMod_{V(L)}$ can be resolved in $(\sfsMod_{V(L)})_{\ev}$ by free $V(L)$-modules, it follows that the forgetful functor induces for each $n \in \N$ a canonical identification
	\begin{equation} \label{eq:Extsmod=Extmod}
	\Ext_{\sfsMod}^n(M,N) = \Ext_{\sfMod}^n(M,N) \quad \text{for all $M,N \in \sfsMod_{V(L)}$.}
	\end{equation}
From now on, we will simply write $\Hom_{V(L)}(M,N)$ and $\Ext_{V(L)}^n(M,N)$ to denote the common values indicated by \eqref{eq:Homsmod=Hommod} and \eqref{eq:Extsmod=Extmod}.

The identification in \eqref{eq:Extsmod=Extmod} is also compatible with cup products. This can be seen using the description of cup products given, for example, in \cite[\S9.3]{Witherspoon:2019}. The key point to check to ensure that this description makes sense for $\sfsMod_{V(L)}$ is that if $P$ and $Q$ are projective modules in $\sfsMod_{V(L)}$, then the tensor product $P \otimes Q$ with the diagonal $V(L)$-action is also projective in $\sfsMod_{V(L)}$. It suffices to assume that $P$ is a free $V(L)$-supermodule, and then the proof of \cite[Lemma 9.2.9]{Witherspoon:2019} produces an isomorphism in $(\sfsMod_{V(L)})_{\ev}$ between $P \otimes Q$ and the free $V(L)$-supermodule $P \otimes Q_{triv}$. Here $Q_{triv}$ denotes the $V(L)$-supermodule having the same underlying superspace as $Q$, but with the action of $V(L)$ defined by the augmentation map $\ve: V(L) \to k$.

By Lemma \ref{lemma:RLSA=RLA}, $V(L)$ identifies with $V(\Lbar)$, the restricted enveloping algebra of the ordinary restricted Lie algebra $\Lbar$, whose underlying Lie algebra is obtained by forgetting the quadratic operator $q: \Lone \to \Lzero$, and whose $p$-operation $z \mapsto z^{\set{2}}$ is defined by \eqref{eq:2-operation}. Then the module categories $\sfMod_{V(L)}$ and $\sfMod_{V(\Lbar)}$ are the same, so by \eqref{eq:Extsmod=Extmod} one gets a canonical identification
	\begin{equation} \label{eq:ExtVL=ExtVLtilde}
	\Ext_{V(L)}^\bullet(M,N) = \Ext_{V(\Lbar)}^\bullet(M,N) \quad \text{for all $M,N \in \sfsMod_{V(L)}$,}
	\end{equation}
which is again evidently compatible with cup products. To emphasize, the left-hand side of \eqref{eq:ExtVL=ExtVLtilde} represents the value calculated in the category of $V(L)$-supermodules, while the right-hand side represents the value calculated without regard for the $\Z_2$-gradings.

\subsection{Support}

Given a Hopf superalgebra $A$ over $k$, the cohomology ring $\Hbul(A,k) = \Ext_A^\bullet(k,k)$ is a graded-commutative superalgebra \cite[Corollary 2.3.6]{Drupieski:2019a}. For $p=2$, this means that $\Hbul(A,k)$ is an ordinary commutative ring. We define the cohomological spectrum of $A$, denoted $\abs{A}$, to be
	\begin{equation} \label{eq:A-spectrum}
	\abs{A} = \Spec \Hbul(A,k),
	\end{equation}
the prime ideal spectrum of $\Hbul(A,k)$, considered as a topological space via the Zariski topology. If the field $k$ is algebraically closed, one can take the maximal ideal spectrum instead, and then $\abs{A}$ will be an affine algebraic variety (assuming that $\Hbul(A,k)$ is a finitely-generated $k$-algebra, which it is for the examples we consider) rather than an affine scheme. We will adopt this latter, more classical perspective later in Section \ref{S:support-for-non-restricted} when discussing support varieties for $U(L)$.

Given $A$-supermodules $M$ and $N$, let $I_A(M,N)$ be the annihilator ideal for the cup product action of $\Hbul(A,k)$ on $\Ext_A^\bullet(M,N)$. The ideal $I_A(M,N)$ is bi-homogeneous (i.e., homogeneous with respect to both the `external' cohomological degree and the `internal' super degree) because the cup product action is compatible with both the cohomological and super degrees. The relative support variety of the pair $(M,N)$ is then defined to be
	\begin{equation} \label{eq:A-relative-support}
	\abs{A}_{(M,N)} = \Spec\left( \Hbul(A,k)/I_A(M,N) \right),
	\end{equation}
a Zariski closed subset of $\abs{A}$, and the support variety of $M$ is defined to be
	\begin{equation} \label{eq:A-support-variety}
	\abs{A}_M = \abs{A}_{(M,M)} = \Spec\left( \Hbul(A,k)/I_A(M,M) \right).
	\end{equation}
For more details on basic notions concerning support varieties for Hopf (super)algebras, see for example \cite[\S2.3]{Drupieski:2019a} or \cite[\S8.3]{Witherspoon:2019}.

Taking $M = N = k$ in \eqref{eq:ExtVL=ExtVLtilde}, the identification of cohomology rings $\Hbul(V(L),k) = \Hbul(V(\Lbar),k)$ implies the equality of cohomological spectra:
	\[
	\abs{V(L)} = \abs{V(\Lbar)}.
	\]
By \cite[Lemma 1.6]{Suslin:1997} and \cite[Theorem 5.2]{Suslin:1997a}, there is a finite universal homeomorphism between $\abs{V(\Lbar)}$ and the restricted nullcone of $\Lbar$:
	\begin{equation} \label{eq:Lbar-spectrum}
	\abs{V(\Lbar)} \simeq \calN_1(\Lbar) := \set{ z \in \Lbar : z^{\set{2}} = 0}.\footnote{The results in \cite{Suslin:1997,Suslin:1997a} are stated in the language of finite group schemes, so the reader must make the standard translation from the category of $V(\Lbar)$-modules to the equivalent category of rational modules for the finite group scheme $G$ whose coordinate algebra $k[G]$ is $V(\Lbar)^*$, the Hopf algebra dual to $V(\Lbar)$.}
	\end{equation}
Writing $z = z_{\zero} + z_{\one}$ for the decomposition of $z$ into its even and odd parts, one has $z^{\set{2}} = 0$ if and only if $(z_{\zero})^{[2]} + q(z_{\one}) + [z_{\one},z_{\zero}] = 0$. Since $(z_{\zero})^{[2]} + q(z_{\one})$ is of even superdegree and $[z_{\one},z_{\zero}]$ is of odd superdegree, this implies that
	\begin{equation} \label{eq:z2=0}
	z^{\set{2}} = 0 \quad \text{if and only if} \quad (z_{\zero})^{[2]} + q(z_{\one}) = 0 \text{ and } [z_{\one},z_{\zero}] = 0.
	\end{equation}
Thus, for the spectrum of $V(L)$ one gets
	\begin{equation} \label{eq:V(L)spectrum}
	\abs{V(L)} \simeq \set{ z = z_{\zero}+z_{\one} \in L : (z_{\zero})^{[2]} + q(z_{\one}) = 0 \text{ and } [z_{\one},z_{\zero}] = 0}.
	\end{equation}
This extends to $p = 2$ the description of $\abs{V(L)}$ we previously obtained for characteristic $p \geq 3$ in \cite[Proposition 5.4.4]{Drupieski:2019a}. In fact, the result here for $p=2$ is stronger than the result in \cite{Drupieski:2019a}, since the identification there is only proved up to a finite morphism of varieties.

Next, the compatibility of \eqref{eq:ExtVL=ExtVLtilde} with cup products implies that
	\begin{equation} \label{eq:L-supp=Lbar-supp}
	\abs{V(L)}_{M} = \abs{V(\Lbar)}_{M} \quad \text{for all $M \in \sfsMod_{V(L)}$.}
	\end{equation}
By \cite[Corollary 6.3.1]{Suslin:1997a}, the homeomorphism \eqref{eq:Lbar-spectrum} restricts for $M$ finite-dimensional to
	\begin{equation} \label{eq:Lbar-supp=rank}
	\abs{V(\Lbar)}_M \simeq \set{z \in \Lbar : z^{\set{2}} = 0 \text{ and } M|_{\subgrp{z}} \text{ is not free}} \cup \set{0}.\footnote{This identification is true as written if the field $k$ is algebraically closed, and if one defines $\abs{V(\Lbar)}$ and support varieties in terms of the maximal ideal spectrum. But if one is working at the level of schemes, the precise statement of the identification involves certain field extensions; see \cite[Proposition 6.1]{Suslin:1997a}. For simplicity of exposition, we have suppressed the extra details.}
	\end{equation}
Here $M|_{\subgrp{z}}$ denotes the restriction of $M$ to the $k$-subalgebra of $V(\Lbar)$ generated by $z$. For $z \neq 0$ this subalgebra has the form $k[z]/(z^2)$, while for $z = 0$ the subalgebra is just the field $k$. (This subalgebra is not a sub-\emph{superalgebra} of $V(\Lbar)$ unless $z$ is homogeneous.) An immediate consequence of \eqref{eq:L-supp=Lbar-supp}, \eqref{eq:Lbar-supp=rank}, and the representation theory of $k[z]/(z^2)$ is the following (cf.\ \cite[Theorem 7.2]{Suslin:1997a}):

\begin{proposition}[Tensor Product Property] \label{prop:tensor-product-property}
Let $L$ be a finite-dimensional restricted Lie super\-algebra over an algebraically closed field of characteristic $p=2$, and let $M$ and $N$ be finite-dimensional $V(L)$-supermodules. Then
	\[
	\abs{V(L)}_{M \otimes N} = \abs{V(L)}_M \cap \abs{V(L)}_N.
	\]
\end{proposition}

Next we give an alternative characterization of the condition that $M|_{\subgrp{z}}$ is not free.
	
\begin{lemma} \label{lemma:notfree-projdiminf}
Let $M$ be a finite-dimensional $V(\Lbar)$-module, and let $z = z_{\zero} + z_{\one} \in \Lbar$ with $z^{\set{2}} = 0$. Set $P = k[u,v]/(u^2+v^2)$, and let $\sigma_z : P \to V(\Lbar)$ be the homomorphism defined by $\sigma_z(u) = z_{\zero}$ and $\sigma_z(v) = z_{\one}$. Write $M\da{\sigma_z}$ for the $P$-module obtained by pulling back $M$ along $\sigma_z$. Then for $z \neq 0$,
	\begin{center}
	$M|_{\subgrp{z}}$ is not free \qquad if and only if \qquad $\projdim_P(M\da{\sigma_z}) = \infty$.
	\end{center}
For $z = 0$, one has that $M|_{\subgrp{z}}$ is free (i.e., $M$ is free over $k$) and $\projdim_P(M\da{\sigma_z}) = \infty$.
\end{lemma}

\begin{proof}
Suppose $z \neq 0$, and let $Q$ be the $k$-subalgebra of $P$ generated by $u+v$. Then $Q \cong k[z]/(z^2)$, so $M|_{\subgrp{z}}$ is not free if and only if $\projdim_Q(M\da{\sigma_z}) = \infty$. Next, $P = Q[u]$ is a polynomial extension of the ring $Q$, so $\flatdim_P(M\da{\sigma_z}) = \infty$ if and only if $\flatdim_Q(M\da{\sigma_z}) = \infty$ by \cite[\S5.3]{Avramov:2018}. Finally, $\projdim_P(M\da{\sigma_z}) = \flatdim_P(M\da{\sigma_z})$ and $\projdim_Q(M\da{\sigma_z}) = \flatdim_Q(M\da{\sigma_z})$ because $P$ and $Q$ are commutative Noetherian rings over which $M$ is finitely generated. If $z = 0$, then $\sigma_z: P \to k$ is the augmentation map, and $\projdim_P(M\da{\sigma_z}) = \projdim_P(k) = \infty$.
\end{proof}

\begin{proposition} \label{prop:rank-projdim}
Let $L$ be a finite-dimensional restricted Lie superalgebra over an algebraically closed field of characteristic $p=2$, and let $M$ be a finite-dimensional $V(L)$-supermodule. Let $P$ and $\sigma_z: P \to V(\Lbar)$ be as defined in Lemma \ref{lemma:notfree-projdiminf}. Then there exists a homeomorphism
	\begin{equation} \label{eq:V(L)modulesupport}
	\abs{V(L)}_M \simeq \set{ z \in L : (z_{\zero})^{[2]} + q(z_{\one}) = 0 \text{ and } [z_{\one},z_{\zero}] = 0 \text{ and } \projdim_P(M\da{\sigma_z}) = \infty}.
	\end{equation}
A $V(L)$-supermodule $M$ is projective if and only if $\abs{V(L)}_M = \set{0}$, or equivalently, if and only if $\projdim_P(M\da{\sigma_z}) < \infty$ for all nonzero $z \in L$ such that $z^{\set{2}} = 0$.
\end{proposition}

\begin{proof}
The asserted homeomorphism follows immediately from \eqref{eq:z2=0}, \eqref{eq:L-supp=Lbar-supp}, \eqref{eq:Lbar-supp=rank}, and the lemma. Next, the fact that the forgetful functor $\Forget: \sfsMod_{V(L)} \to \sfMod_{V(L)} = \sfMod_{V(\Lbar)}$ is surjective on morphisms implies that if a $V(L)$-supermodule is projective as an object in $\sfMod_{V(\Lbar)}$, then it is also projective as an object in (the underlying even subcategory of) $\sfsMod_{V(L)}$. Then the last sentence of of the proposition is true by the well-known fact that a finite-dimensional $V(\Lbar)$-module is projective if and and only if $\abs{V(\Lbar)}_M = \set{0}$; cf.\ \cite[Proposition 1.5]{Friedlander:1987}.
\end{proof}

\begin{remark}\label{rem:detecting-finite=proj-dim}
In \cite[Lemma 3.1.3]{Drupieski:2021b}, we interpreted our previous work from  \cite{Drupieski:2021a} to show that an identification analogous to \eqref{eq:V(L)modulesupport} holds for the support varieties of arbitrary finite-dimensional \emph{$p$-nilpotent} restricted Lie super\-algebras over (algebraically closed) fields of characteristic $p \geq 3$. In \cite[\S1.2]{Drupieski:2021a}, we conjectured that a similar identification should hold for all finite-dimensional restricted Lie superalgebras (and more generally, for all infinitesimal supergroup schemes) in characteristic $p \geq 3$. The primary obstruction to proving such a generalization in characteristic $p \geq 3$ is the lack of a characterization of projectivity similar to that stated in the last sentence of Proposition \ref{prop:rank-projdim}.
\end{remark}

Recall that a closed set $W \subseteq \abs{A}$ is \emph{conical} if it is defined by a homogeneous (respect to the cohomological degree) ideal in $\Hbul(A,k)$. Given a Hopf superalgebra $A$, let us say that a closed set $W \subseteq \abs{A}$ is \emph{super-conical} if it is defined by a bi-homogeneous ideal in $\Hbul(A,k)$, i.e., an ideal that is homogeneous with respect to both the `external' cohomological degree and the `internal' super degree. Given a $V(L)$-supermodule $M$, the support variety $\abs{V(L)}_M$ is super-conical because its defining ideal $I_{V(L)}(M,M)$ is automatically bi-homogeneous. The next proposition states that every super-conical subset of $\abs{V(L)}$ is of this form.

\begin{proposition}[Realizability] \label{prop:realizability}
Let $L$ be a finite-dimensional restricted Lie superalgebra over an algebraically closed field of characteristic $2$. Then $W \subseteq \abs{V(L)}$ is the support variety $\abs{V(L)}_M$ of a finite-dimensional $V(L)$-supermodule if and only if $W$ is a closed, super-conical subset of $\abs{V(L)}$.
\end{proposition}

\begin{proof}
We have already explained the `only if' portion of the statement, so we just need to prove the `if' portion. Let $W \subseteq \abs{V(L)}$ be a closed, super-conical subset. Then $W$ is defined by a bi-homo\-geneous ideal in $\Hbul(V(L),k)$, which in turn is generated by a finite set of bi-homogeneous elements. By the Tensor Product Property (Proposition \ref{prop:tensor-product-property}), it suffices to show that if $W = Z(\zeta)$ is the closed set defined by a single bi-homogeneous element $\zeta \in \opH^n(V(L),k)$, then $W = \abs{V(L)}_M$ for some finite-dimensional $V(L)$-supermodule $M$.

Choose a projective resolution $(P_\bullet,d)$ in $(\sfsMod_{V(L)})_{\ev}$ of the trivial $V(L)$-supermodule $k$, and a homogeneous (with respect to the super degree) cocycle representative $g: P_n \to k$ for $\zeta$. Let $M$ be the kernel of the induced map $\wh{g}: P_n/d(P_{n+1}) \to k$. (As an ungraded module in $\sfMod_{V(\Lbar)}$, and perhaps modulo a projective summand, $P_n/d(P_{n+1})$ is the $n$-th syzygy module $\Omega^n(k)$ for the trivial $V(\Lbar)$-module $k$.) Then $M$ is a $V(L)$-supermodule because $d: P_{n+1} \to P_n$ and $g: P_n \to k$ are both homogeneous $V(L)$-super\-module homomorphisms. Now the proof of \cite[Theorem 7.5]{Suslin:1997a} (cf.\ also \cite[\S4]{Friedlander:1987}) shows that $\abs{V(\Lbar)}_M = Z(\zeta)$, hence $\abs{V(L)}_M = Z(\zeta)$ by \eqref{eq:L-supp=Lbar-supp}.
\end{proof}

\subsection{\texorpdfstring{$\Proj(A)$ and $\Proj_s(A)$}{Proj(A) and Projs(A)}} \label{subsec:proj-and-projs}

In anticipation of Section \ref{subsec:tensor-triangular-geometry}, we now discuss a modification of the projective spectrum that is useful for graded(-commutative) superalgebras in characteristic $2$.  If we were in characteristic other than $2$, then $\Z_{2}$-gradings could be encoded by an action of the cyclic group $C_{2}$ and the following would essentially be the theory of $C_{2}\text{-}\Proj (A)$ as described in the literature; see for example \cite{Lorenz:2009}. Since we are in characteristic $2$, and for the sake of completeness, we provide the relevant details.

Let $A$ be a commutative $k$-algebra that is graded by $\N \times \Z_2$. Let $A = \bigoplus_{i \in \N} A^i$ be the decomposition of $A$ into its $\N$-graded components, and assume that $A^0 = k$. Let $\Proj(A)$ be the subset of $\Spec(A)$ consisting of the prime ideals that are homogeneous with respect to the $\N$-grading on $A$, and which are properly contained in the maximal ideal $A^+ = \bigoplus_{i \geq 1} A^i$. A subset $W$ of $\Proj(A)$ is closed in the Zariski topology if there exists an $\N$-homogeneous ideal $I$ of $A$ such that
	\[
	W = Z(I) := \set{ \fp \in \Proj(A) : I \subseteq \fp}.
	\]

We say that an ideal $I$ of $A$ is \emph{bi-homogeneous} if it is homogeneous with respect to the full $\N \times \Z_2$-grading on $A$. For example, the ideal $A^+$ is bi-homogeneous. Given an ideal $I$ of $A$, let $I_s \subseteq I$ be the largest bi-homogeneous subideal of $I$. Equivalently, $I_s$ is the ideal in $A$ generated by the bi-homogeneous elements of $I$. Of course, $I = I_s$ if $I$ is bi-homogeneous. If $J$ is a bi-homogeneous ideal and $I$ is an arbitrary ideal, then $J \subseteq I$ if and only if $J \subseteq I_s$. If $\set{ I_\lambda : \lambda \in \Omega}$ is any family of ideals in $A$, then $(\bigcap_{\lambda \in \Omega} I_\lambda)_s = \bigcap_{\lambda \in \Omega} (I_\lambda)_s$.

Given an ideal $I$ of $A$, let $\sqrt{I}$ be the usual radical of $I$, and set $\sqrt[s]{I} = (\sqrt{I})_s$. Then $I_s \subseteq \sqrt[s]{I}$, and $\sqrt[s]{I}$ has the property that if $a \in A$ is a bi-homogeneous element with $a^n \in I_s$ for some $n \geq 1$, then $a \in \sqrt[s]{I}$. We say that a bi-homogeneous ideal $J$ is an \emph{$s$-radical ideal} if $J = \sqrt[s]{J}$. We say that a bi-homogeneous ideal $P$ is \emph{$s$-prime} if whenever $I$ and $J$ are bi-homogeneous ideals with $IJ \subseteq P$, then either $I \subseteq P$ or $J \subseteq P$. Equivalently, $P$ is $s$-prime if whenever $a$ and $b$ are bi-homogeneous elements with $ab \in P$, either $a \in P$ or $b \in P$. If $P$ is a prime ideal in $A$, then $P_s$ is an $s$-prime ideal. Every $s$-prime ideal in $A$ is an $s$-radical ideal.

\begin{lemma}
Let $P$ be an $s$-prime ideal of $A$. Then $\sqrt{P}$ is an $\N$-homogeneous prime ideal of $A$.
\end{lemma}

\begin{proof}
Ignoring the $\Z_2$-gradings, $P$ is an $\N$-homogeneous ideal of $A$, so its radical is also an $\N$-homo\-geneous ideal. To show that $\sqrt{P}$ is prime, let $x,y \in A$, and suppose $xy \in \sqrt{P}$, say with $x^n y^n = (xy)^n \in P$. Replacing $x$ and $y$ with $x^n$ and $y^n$, we may assume that $xy \in P$. Then also $x^2 y^2 = (xy)^2 \in P$. Write $x = \sum x_i$ with $x_i \in A^i$, and write $x_i = x_{i,\zero} + x_{i,\one}$ with $x_{i,\zero} \in (A^i)_{\zero}$ and $x_{i,\one} \in (A^i)_{\one}$. Since $A$ is commutative, the binomial theorem modulo $2$ implies that $x^2 = \sum (x_i)^2 = \sum [ (x_{i,\zero})^2 + (x_{i,\one})^2 ]$. Then $x^2$ is a sum of terms of even superdegree, and we can rewrite $x^2$ as $x^2 = \sum x_{2i}'$, with $x_{2i}' \in (A^{2i})_{\zero}$. Similarly, $y^2 = \sum y_{2i}'$ with $y_{2i}' \in (A^{2i})_{\zero}$.

Now suppose that $x^2 \notin P$ and $y^2 \notin P$. Then there exist minimal indices $s$ and $t$ such that $x_{2s}' \notin P$ and $y_{2t}' \notin P$. Since $x^2$ and $y^2$ are each of even superdegree, the $\N$-graded components of $x^2y^2$ are each bi-homogeneous. Since $P$ is bi-homogeneous, and since $x^2 y^2 \in P$, this means in particular that the component of $\N$-degree $s+t$ in $x^2y^2$ is an element of $P$. This component is equal to $\sum_{\ell \in \Z} x_{2s-\ell}' y_{2t+\ell}'$. By the minimality of $s$ and $t$, each summand $x_{2s-\ell}' y_{2t+\ell}'$ with $\ell \neq 0$ is an element of $P$, because at least one of the factors is an element of $P$ by assumption. Then $x_{2s}' y_{2t}'$ must also be an element of $P$. The factors $x_{2s}'$ and $y_{2t}'$ are each bi-homogeneous, so the fact that $P$ is an $s$-prime ideal implies that either $x_{2s}' \in P$ or $y_{2t}' \in P$, in contradiction to the assumptions on $s$ and $t$. Thus, it must be the case that either $x^2 \in P$ or $y^2 \in P$, and hence $x \in \sqrt{P}$ or $y \in \sqrt{P}$.
\end{proof}

\begin{definition}
Define $\Proj_s(A)$ to be the set of all (bi-homogeneous) $s$-prime ideals in $A$ that are properly contained in the (bi-homogeneous, $s$-prime) maximal ideal $A^+$. We say that a subset $W$ of $\Proj_s(A)$ is \emph{closed} if there exists a bi-homogeneous ideal $I$ of $A$ such that
	\[
	W = Z_s(I) := \set{ P \in \Proj_s(A) : I \subseteq P}.
	\]
\end{definition}

The properties listed below are now straightforward to verify. Unless stated otherwise, $I$ and $J$ denote bi-homogeneous ideals in $A$.
	\begin{enumerate} [label={(P\arabic*)}]
	\item \label{item:unions} $Z_s(I) \cup Z_s(J) = Z_s(IJ)$.
	\item \label{item:intersections} $Z_s(I) \cap Z_s(J) = Z_s(I+J)$.
	\item \label{item:P-to-Ps} If $\fp \in \Proj(A)$, then $\fp_s \in \Proj_s(A)$.
	\item \label{item:continuous} If $\varphi: \Proj(A) \to \Proj_s(A)$ is the function $\varphi(\fp) = \fp_s$, then $\varphi^{-1}(Z_s(I)) = Z(I)$.
	\item If $P \in \Proj_s(A)$, then $\sqrt{P} \in \Proj(A)$.
	\item \label{item:section} If $\psi: \Proj_s(A) \to \Proj(A)$ is the function $\psi(P) = \sqrt{P}$, then $(\varphi \circ \psi)(P) = \sqrt[s]{P} = P$.
	\end{enumerate}
The first two properties imply that finite unions and arbitrary intersections of closed sets are again closed. Thus the closed subsets of $\Proj_s(A)$ define a topology, which we call the Zariski topology on $\Proj_s(A)$. Property \ref{item:continuous} then implies that $\varphi: \Proj(A) \to \Proj_s(A)$ is a continuous function.

If $K$ is an arbitrary $\N$-homogeneous ideal in $A$, one has $\varphi(Z(K)) \subseteq Z_s(K_s)$. On the other hand, if $I$ is a bi-homogeneous ideal, and if $P \in Z_s(I)$, then $I \subseteq P \subseteq \sqrt{P}$. Then property \ref{item:section} implies:
	\begin{enumerate}[resume,label={(P\arabic*)}]
	\item \label{item:bi-homog-closed-to-closed} If $I$ is bi-homogeneous, then $\varphi(Z(I)) = Z_s(I)$.
	\end{enumerate}
Next observe that if $Z_s(I) \subseteq Z_s(J)$, then
	\[
	Z(I) = \varphi^{-1}(Z_s(I)) \subseteq \varphi^{-1}(Z_s(J)) = Z(J).
	\]
But $Z(I) \subseteq Z(J)$ implies that $\sqrt{I} \supseteq \sqrt{J}$, and hence $\sqrt[s]{I} \supseteq \sqrt[s]{J}$. Now one can check that:
	\begin{enumerate}[resume,label={(P\arabic*)}]
	\item \label{item:define-by-radical} $Z_s(I) = Z_s(\sqrt[s]{I})$.
	\item \label{item:ideal-to-Z} If $I \subseteq J$, then $Z_s(I) \supseteq Z_s(J)$.
	\item \label{item:Z-to-radical} If $Z_s(I) \subseteq Z_s(J)$, then $\sqrt[s]{I} \supseteq \sqrt[s]{J}$.
	\item \label{item:noetherian} If $A$ is a Noetherian ring, then $\Proj_s(A)$ is a Noetherian topological space, i.e., if $Z_s(I_1) \supseteq Z_s(I_2) \supseteq Z_s(I_3) \supseteq \cdots$ is a descending chain of closed sets, then there exists an integer $m$ such that $Z_s(I_j) = Z_s(I_m)$ for all $j \geq m$. 
	\end{enumerate}
	
\begin{lemma} \label{lemma:irreducible-closure-s-prime}
If $W = Z_s(I)$ is an irreducible closed set in $\Proj_s(A)$, then $\sqrt[s]{I}$ is an $s$-prime ideal. In particular, if $W$ is nonempty, then $W$ is the closure of the point $\sqrt[s]{I} \in \Proj_s(A)$, and $\sqrt[s]{I}$ is the unique point in $\Proj_s(A)$ with this property.
\end{lemma}

\begin{proof}
To show that $\sqrt[s]{I}$ is an $s$-prime ideal, one can apply the usual argument for $\Spec(A)$; see for example the proof of \cite[Proposition 4.5.3]{Cox:2015}. Then $W = Z_s(I) = Z_s(\sqrt[s]{I})$ by property \ref{item:define-by-radical} above, while property \ref{item:Z-to-radical} implies that if $W = \overline{\set{P}} = Z_s(P)$ for some $P \in \Proj_s(A)$, then $P = \sqrt[s]{I}$.
\end{proof}

Combining property \ref{item:noetherian} and Lemma \ref{lemma:irreducible-closure-s-prime}, we get:
	\begin{enumerate}[resume,label={(P\arabic*)}]
	\item If $A$ is a Noetherian ring, then $\Proj_s(A)$ is a Zariski space, i.e., it is a Noetherian topological space, and every nonempty closed irreducible subset is the closure of a unique point.
	\end{enumerate}

\subsection{Tensor triangular geometry} \label{subsec:tensor-triangular-geometry}

Now take $A = \Hbul(V(L),k)$. In this section we discuss how the geometry of $\Proj_s(A) = \Proj_s (\Hbul(V(L),k))$ classifies certain categories of $V(L)$-supermodules. First we recall the corresponding results for the category of ungraded $V(\Lbar)$-modules.

Let $\mathbf{K}$ be a symmetric tensor triangulated category. Recall that a triangulated subcategory of $\mathbf{K}$ is \emph{thick} if it is closed under taking direct summands, and is a \emph{tensor ideal} (or \emph{$\otimes$-ideal}) if it is closed under $M \otimes -$ for any object $M$ in $\mathbf{K}$. A $\otimes$-ideal $\calJ$ is \emph{radical} if whenever $M^{\otimes n} \in \calJ$ for some $n \geq 1$, then also $M \in \calJ$. A thick $\otimes$-ideal $\calP$ is \emph{prime} if whenever $M,N \in \mathbf{K}$ and $M \otimes N \in \calP$, then either $M \in \calP$ or $N \in \calP$.  The \emph{Balmer spectrum} of $\mathbf{K}$, denoted $\Spc (\mathbf{K})$, is the set of all proper prime thick $\otimes$-ideals of $\mathbf{K}$, considered as a topological space via the Zariski topology. Write $2^{\Spc (\mathbf{K})}$ for the power set of $\Spc (\mathbf{K})$. Then the support function 
	\[
	\supp ' : \mathbf{K} \to 2^{\Spc (\mathbf{K})}
	\]
is defined by $\supp' (M)  = \set{\calP \in \Spc (\mathbf{K}) : M \notin \calP }$ for each object $M \in \mathbf{K}$.

Let $\stmodVLbar$ (resp., $\StModVLbar$) be the stable module category of finite-dimensional (resp., all) $V(\Lbar)$-modules. These are the categories obtained from $\sfmod_{V(\Lbar)}$ and $\sfMod_{V(\Lbar)}$, respectively, by quotienting out the morphisms that factor through a projective. The objects of the stable category are the same as in the corresponding ordinary module category, but with projective objects now isomorphic to zero. The categories $\stmodVLbar$ and $\StModVLbar$ are tensor triangulated categories, with suspension $\Sigma$ equal to $\Omega^{-1}(-)$, the inverse of the syzygy functor. The category $\StModVLbar$ is compactly generated, with subcategory of compact objects $(\StModVLbar)^c = \stmodVLbar$.

Now set $X = \Proj (\Hbul(V(\Lbar),k))$. In \cite{Benson:2018}, Benson, Iyengar, Krause, and Pevtsova use localization functors to define a support function
	\begin{equation} \label{eq:supp}
	\supp: \StModVLbar \to 2^X
	\end{equation}
that satisfies the properties listed below.\footnote{Again, the results in \cite{Benson:2018} are stated in the language of finite group schemes, so the reader must translate to the category of rational $G$-modules, where $k[G] = V(\g)^*$.} Except when indicated, $M$, $M_i$, and $N$ denote arbitrary objects in $\StModVLbar$. 
	\begin{enumerate} [label={(S\arabic*)}] \setcounter{enumi}{-1}
	\item \label{item:support-zero} $\supp(0) = \emptyset$

	\item $\supp(k) = X$

	\item $\supp(M \oplus N) = \supp(M) \cup \supp(N)$

	\item $\supp(\Sigma M) = \supp(M)$

	\item For each distinguished triangle $M_1 \to M_2 \to M_3 \to \Sigma M_1$ in $\StModVLbar$,
		\[
		\supp(M_2) \subseteq \supp(M_1) \cup \supp(M_3)
		\]

	\item \label{item:tensor-product} Tensor Product Property: For all compact $M$ and arbitrary $N$,
		\[
		\supp(M \otimes N) = \supp(M) \cap \supp(N).\footnote{Actually, this property holds without the assumption that $M$ is compact, but we will find it convenient to use this more restricted formulation of the Tensor Product Property.}
		\]

	\item $\supp(\bigoplus_{i \in I} M_i) = \bigcup_{i \in I} \supp(M_i)$ for any set $I$.

	\item \label{item:compact-has-closed-support} For all compact $M$, $\supp(M)$ is Zariski closed with quasi-compact complement.\footnote{As used in algebraic geometry, the term \emph{quasi-compact} means that every open cover has a finite subcover.}

	\item \label{item:compact-realization} For every Zariski closed subset $W \subseteq X$ with quasi-compact complement, there exists a compact object $M$ such that $\supp(M) = W$.

	\item \label{item:support-detect-projectivity} If $\supp(M) = \emptyset$, then $M \cong 0$.

	\item For each finite-dimensional $V(\Lbar)$-module $M$, $\supp(M^*) = \supp(M)$.
	\end{enumerate}
We note that the `quasi-compact complement' condition in \ref{item:compact-has-closed-support} and \ref{item:compact-realization} is automatically satisfied by any Zariski closed subset of $X$, as a consequence of the fact that $\Hbul(V(\Lbar),k)$ is a finitely-generated commutative $k$-algebra and hence a Noetherian ring. Elaborating further on \ref{item:compact-has-closed-support}, if $M$ is a compact object in $\StModVLbar$ (equivalently, if $M$ is isomorphic, modulo projective summands, to a finite-dimensional $V(\Lbar)$-module), then $\supp(M)$ is the Zariski closed subset of $X$ defined by $I(M)$, the annihilator ideal for the cup product action of $\Hbul(V(\Lbar),k)$ on $\Ext_{V(\Lbar)}^\bullet(M,M)$.

Properties \ref{item:support-zero}--\ref{item:support-detect-projectivity} imply by Theorem 1.5 of \cite{DellAmbrogio:2010} that the restriction of \eqref{eq:supp} to $\stmodVLbar$ is a classifying support datum in the sense of \cite[Theorem 2.19]{DellAmbrogio:2010}; cf.\ also Theorems 3.4.1 and 3.5.1 of \cite{Boe:2017}. Then there is a canonical homeomorphism from $X = \Proj (\Hbul(V(\Lbar),k))$ to the Balmer spectrum of the category $\stmodVLbar$,
	\begin{equation} \label{eq:Balmer-spec-homeo}
	\Proj (\Hbul(V(\Lbar),k)) \simeq \Spc( \stmodVLbar ),
	\end{equation}
which identifies the support functions $\supp$ and $\supp'$. Furthermore, there are inverse bijections
	\begin{multline} \label{eq:BIKP-Theorem10.3}
	\set{ \text{specialization closed subsets $\calV$ of $\Proj (\Hbul(V(\Lbar),k))$}} \overset{\Gamma}{\underset{\Theta} \rightleftarrows} \\
	\set{ \text{thick $\otimes$-ideal subcategories $\calJ$ of $\stmodVLbar$}},
	\end{multline}
defined by $\Gamma(\calV) = \set{ M \in \stmodVLbar : \supp(M) \subseteq \calV}$ and $\Theta(\calJ) = \bigcup_{M \in \calJ} \supp(M)$. (Recall that a subset of a topological space is \emph{specialization closed} if it is a union of arbitrarily many closed sets.) We note that \cite[Theorem 2.19]{DellAmbrogio:2010} is stated in terms of \emph{radical} thick $\otimes$-ideals, but using the fact that each finite-dimensional $V(\Lbar)$-module $M$ is a direct summand of $M \otimes M^* \otimes M$, one can see that every thick $\otimes$-ideal of $\stmodVLbar$ is automatically radical.

Now we consider the super setting. Recall from Section \ref{subsec:restricted-cohomology} that $(\sfsMod_{V(L)})_{\ev}$ is a Frobenius category. Then we can consider the corresponding stable category (see, e.g., \cite[I.2]{Happel:1988}), which we denote $\StsModVL$. It is a compactly generated tensor-triangulated category, with subcategory of compact objects $(\StsModVL)^c = \stsmodVL$, the stable category of $(\sfsmod_{V(L)})_{\ev}$. The forgetful functor $\sfsMod_{V(L)} \to \sfMod_{V(L)}$ induces a (exact) triangulated functor
	\[
	\Forget: \StsModVL \to \StModVL = \StModVLbar.
	\]
In particular, the forgetful functor maps projective objects in $\sfsMod_{V(L)}$ to projective objects in $\sfMod_{V(L)}$. Via this functor, we will consider objects from $\StsModVL$ as objects in $\StModVL$.

Our goal now is to establish analogues of \eqref{eq:Balmer-spec-homeo} and \eqref{eq:BIKP-Theorem10.3} for $\stsmod_{V(L)}$. Set
	\[
	X_s = \Proj_s (\Hbul(V(L),k)),
	\]
and let 
\[
\varphi:  \Proj (\Hbul(V(L),k)) \to \Proj_s (\Hbul(V(L),k))
\] 
be the function $\fp \mapsto \fp_s$ as defined in \ref{item:continuous} above. Using the fact that $\Hbul(V(L),k)$ is a Noetherian ring and the properties listed in Section \ref{subsec:proj-and-projs}, one can again deduce that every Zariski closed set in $X_s$ has quasi-compact complement. Define 
	\[
	\supp_s : \StsMod_{V(L)} \to 2^{X_s}
	\]
by
	\[
	\supp_s(M) = \varphi\left( \supp(\Forget (M)) \right).
	\]

\begin{theorem} \label{T:BalmerSpectrum}
Let $L$ be a finite-dimensional restricted Lie superalgebra over a field $k$ of characteristic $p = 2$. There is a canonical homeomorphism from $X_s = \Proj_s (\Hbul(V(L),k))$ to the Balmer spectrum of the category $\stsmodVL$,
	\begin{equation} \label{eq:Balmer-spec-homeo-s}
	\Proj_s (\Hbul(V(\Lbar),k)) \simeq \Spc( \stmodVL ),
	\end{equation}
which identifies $\supp_{s}$ with $\supp'$. Furthermore, there are inverse bijections
	\begin{multline} \label{eq:BIKP-Theorem10.3-s}
	\set{ \text{specializataion closed subsets $\calV$ of $\Proj_s(\Hbul(V(L),k))$}} \overset{\Gamma}{\underset{\Theta} \rightleftarrows} \\
	\set{ \text{thick $\otimes$-ideal subcategories $\calJ$ of $\stsmodVL$}},
	\end{multline}
defined by $\Gamma(\calV) = \set{ M \in \stsmodVL : \supp_s(M) \subseteq \calV}$ and $\Theta(\calJ) = \bigcup_{M \in \calJ} \supp_s(M)$.
\end{theorem}

\begin{proof}
The theorem follows directly from \cite[Theorem 1.5]{DellAmbrogio:2010}, once we verify that the analogues of \ref{item:support-zero}--\ref{item:support-detect-projectivity} hold for $\supp_s$. Many of the properties are immediate consequences of the fact that the forgetful functor $\Forget: \StsModVL \to \StModVL$ is a triangulated functor, and the fact that the function $\varphi: X \to X_s$ is a continuous surjection. We will make a few comments about those properties that are not entirely immediate.

\ref{item:compact-has-closed-support}: Let $M$ be a finite-dimensional $V(L)$-supermodule. Then $\supp(M) = Z(I(M))$, the Zariski closed subset of $\Proj \Hbul(V(L),k)$ defined by $I(M)$, the annihilator ideal for the cup product action of $\Hbul(V(L),k)$ on $\Ext_{V(L)}^\bullet(M,M)$. The ideal $I(M)$ is bi-homogeneous in $\Hbul(V(L),k)$, so then
	\[
	\supp_s(M) = \varphi\left( Z(I(M)) \right) = Z_s(I(M))
	\]
by property \ref{item:bi-homog-closed-to-closed}. So $\supp_s(M)$ is Zariski closed (with quasi-compact complement) in $X_s$.

\ref{item:tensor-product}: Let $M$ be a finite-dimensional $V(L)$-supermodule, and let $N$ be an arbitrary $V(L)$-super\-module. By property \ref{item:tensor-product} for $\supp$, one has $\supp(M \otimes N) = \supp(M) \cap \supp(N)$. Then
	\[
	\supp_s(M \otimes N) = \varphi\left( \supp(M) \cap \supp(N) \right) \subseteq \supp_s(M) \cap \supp_s(N).
	\]
For the reverse set inclusion, let $P \in \supp_s(M) \cap \supp_s(N)$. Then there exist $\fp \in \supp(M)$ and $\fq \in \supp(N)$ such that $\varphi(\fp) = P = \varphi(\fq)$, i.e., $\fp_s = P = \fq_s$. Since $\supp(M) = Z(I(M))$ and $I(M)$ is bi-homogeneous, this implies that $I(M) \subseteq \fp_s = \fq_s \subseteq \fq$. Then $\fq \in \supp(M) \cap \supp(N) = \supp(M \otimes N)$, and $P = \fq_s = \varphi(\fq) \in \supp_s(M \otimes N)$.

\ref{item:compact-realization}: Let $W \subseteq X_s$ be a Zarsiki closed set (with quasi-compact complement). Then $W = Z_s(I)$ for some bi-homogeneous ideal $I$ in $\Hbul(V(L),k)$. By Proposition \ref{prop:realizability}, there exists a finite-dimensional $V(L)$-supermodule $M$ such that $\supp(M) = Z(I)$. Then $\supp_s(M) = Z_s(I) = W$ by \ref{item:bi-homog-closed-to-closed}.

\ref{item:support-detect-projectivity}: If $M \cong 0$ in $\StsModVL$, then $\supp_s(M) = \emptyset$ by property \ref{item:support-zero} for $\supp$ and the fact that the forgetful functor maps the zero object of $\StsModVL$ to the zero object of $\StModVL$. Conversely, if $\supp_s(M) = \emptyset$, then $\supp(M) = \emptyset$, and hence $M \cong 0$ in $\StModVL$ by property \ref{item:support-detect-projectivity} for $\supp$. Then $M$ is projective as an object in $\sfMod_{V(L)}$. Since the forgetful functor $\sfsMod_{V(L)} \to \sfMod_{V(L)}$ is exact and fully faithful, this implies that $M$ is projective in $(\sfsMod_{V(L)})_{\ev}$, and hence $M \cong 0$ in $\StsModVL$. Thus, $\supp_s(M) = \emptyset$ if and only if $M \cong 0$ in $\StsModVL$.
\end{proof}

\section{Support for non-restricted Lie superalgebras in characteristic 2} \label{S:support-for-non-restricted}

\subsection{Lie superalgebra cohomology} 

Given a $k$-vector space $V$, let $S(V)$, $\Lambda(V)$, and $\Gamma(V)$ be the symmetric, exterior\footnote{For $p=2$, one has $\Lambda(V) \cong S(V)/\subgrp{x^2 : x \in V}$; i.e., $\Lambda(V)$ is isomorphic to a truncated polynomial ring.}, and divided power algebras on $V$, respectively. We denote a typical monomial in $S(V)$ by $v_1^{a_1} \cdots v_t^{a_t}$, and a typical monomial in $\Lambda(V)$ by $\subgrp{v_1 \cdots v_n}$. The divided power algebra $\Gamma(V)$ is generated by the symbols $\gamma_r(v)$ for $r \in \N$ and $v \in V$. These generators satisfy the relations $\gamma_0(v) = 1$, $\gamma_r(0) = 0$ if $r \geq 1$, $\gamma_r(v)\gamma_s(v) = \binom{r+s}{r} \gamma_{r+s}(v)$, and $\gamma_r(v+v') = \sum_{i=0}^r \gamma_i(v)\gamma_{r-i}(v')$. Since $p = \chr(k) = 2$, the algebras $S(V)$, $\Lambda(V)$, and $\Gamma(V)$ are each naturally Hopf algebras: the subspace $V$ consists of primitive elements in each of $S(V)$ and $\Lambda(V)$, and the coproduct on $\Gamma(V)$ satisfies $\Delta(\gamma_r(v)) = \sum_{i+j=r} \gamma_i(v) \otimes \gamma_j(v)$.\footnote{If $p \neq 2$, then $\Lambda(V)$ is not an ordinary Hopf algebra but is instead a Hopf superalgebra.} The algebras $S(V)$, $\Lambda(V)$, and $\Gamma(V)$ are each naturally $\N$-graded, and these $\N$-gradings are compatible with the Hopf algebra structures.

Now let $L = \Lzero \oplus \Lone$ be Lie super\-algebra over $k$, and let $\Ybar(L) = \Lambda(\Lzero) \otimes \Gamma(\Lone)$ be the tensor product of algebras.\footnote{If $p \neq 2$, take the graded tensor product of superalgebras, in which $(a \otimes b)(c \otimes d) = (-1)^{\ol{b} \cdot \ol{c} + \deg(b) \cdot \deg(c)} ac \otimes bd$.} Then $\Ybar(L)$ inherits the structure of $\N$-graded Hopf algebra, with $\Ybar_\ell(L) = \bigoplus_{i+j=\ell} \Lambda^i(\Lzero) \otimes \Gamma^j(\Lone)$. In particular, $\Ybar_1(L) \cong L$, and we write $s: L \rightarrow \Ybar_1(L)$ for the natural identification. The right adjoint action of $L$ on itself extends to a right action of $L$ on $\Ybar(L)$ such that $(z_1z_2).u = z_1(z_2.u) + (-1)^{\ol{u} \cdot \ol{z}_2} (z_1.u)z_2$ for all $z_1,z_2 \in \Ybar(L)$ and $u \in L$ (i.e., $u$ acts by right superderivations), and such that $\subgrp{x}.u = s([x,u])$ and $\gamma_r(y).u = \gamma_{r-1}(y)s([y,u])$ for all $x \in \Lzero$, $y \in \Lone$, and $r \in \N$. The right action of $L$ on $\Ybar(L)$ then extends to a right action of the universal enveloping super\-algebra $U(L)$ on $\Ybar(L)$. Now define $Y(L)$ to be the corresponding smash product algebra $U(L) \# \Ybar(L)$. The coproducts on $U(L)$ and $\Ybar(L)$ induce a bialgebra structure on $Y(L)$. Then $Y(L)$ is an $\N$-graded bialgebra with $U(L)$ concentrated in $\N$-degree $0$. As a vector space and as a left $U(L)$-module, $Y(L) = U(L) \otimes \Ybar(L)$. From now on we denote the product in $Y(L)$ by juxtaposition. Then $Y(L)$ is spanned by the monomials
\begin{equation} \label{eq:standardmonomials}
u\subgrp{x_1 \cdots x_s} \gamma_{a_1}(y_1) \cdots \gamma_{a_t}(y_t),
\end{equation}
with $u \in U(L)$, $x_i \in \Lzero$, $y_j \in \Lone$, and $a_j \in \N$.

\begin{theorem} \label{theorem:Koszulresolution}
Let $L$ be a finite-dimensional Lie superalgebra over the field $k$. There exists a differential $d: Y(L) \rightarrow Y(L)$ (which we call the Koszul differential) making $Y(L)$ into a differential graded bialgebra and into a left $U(L)$-free resolution of the trivial module $k$ (which we call the Koszul resolution of $k$). Given $u \in U(L)$, $x \in \Lzero$, and $y \in \Lone$, one has
\begin{align*}
d(u) &= 0,\\
d(\subgrp{x}) &= x, \quad \text{and} \\
d(\gamma_r(y)) &= y\gamma_{r-1}(y) - \subgrp{q(y)}\gamma_{r-2}(y).
\end{align*}
\end{theorem}

\begin{proof}
Apply the same line of reasoning as in the proof of \cite[Theorem 3.1.1]{Drupieski:2013c}.
\end{proof}

Let $M$ and $N$ be left $U(L)$-supermodules. Given $n \in \N$, set
	\[
	C^n(L,M,N) = \Hom_{U(L)}(Y_n(L) \otimes M,N).
	\]
Here we consider $Y(L) \otimes M$ as a $U(L)$-module via the coproduct, and as a chain complex with differential induced by the Koszul differential $d$ on $Y(L)$. In particular, set $\Cbul(L,k) = \Cbul(L,k,k)$. The Koszul differential induces a differential $\partial$ on $\Cbul(L,M,N)$, making $\Cbul(L,M,N)$ into a cochain complex. Then the Lie superalgebra cohomology group $\Ext_L^\bullet(M,N) = \Ext_{U(L)}^\bullet(M,N)$ can be computed as the $n$-th cohomology group of the complex $\Cbul(L,M,N)$. The coalgebra structure of $Y(L)$ induces the structure of a differential graded algebra on $\Cbul(L,k)$, and the structure of a differential graded $\Cbul(L,k)$-bimodule on $\Cbul(L,M,N)$. Specifically,
	\begin{equation} \label{eq:Koszulcochain}
	\Cbul(L,k) \cong (\Ybar_\bullet(L))^* \cong \Lambda(\Lzero^*) \otimes S(\Lone^*)
	\end{equation}
as an $\N$-graded algebra.\footnote{If $p \neq 2$, this is the graded tensor product of algebras.} Passing to cohomology, one gets the cup product structure on
	\[
	\Hbul(L,k) = \Ext_L^\bullet(k,k) = \Ext_{U(L)}^\bullet(k,k),
	\]
and the $\Hbul(L,k)$-bimodule structure on
	\[
	\Ext_L^\bullet(M,N) = \Ext_{U(L)}^\bullet(M,N).
	\]

As in Section \ref{S:support-for-restricted}, let $\sfsMod_{U(L)}$ be the category of left $U(L)$-supermodules, and let $\sfMod_{U(L)}$ be the category of arbitrary (not necessarily $\Z_2$-graded) left $U(L)$-modules. For each $M \in \sfsMod_{U(L)}$, the tensor product $Y(L) \otimes M$ provides a resolution of $M$ in $(\sfsMod_{U(L)})_{\ev}$ by free $U(L)$-supermodules, so by the same reasoning as for \eqref{eq:Extsmod=Extmod}, one gets
	\begin{equation} \label{eq:U(L)Extsmod=Extmod}
	\Ext_{\sfsMod_{U(L)}}^\bullet(M,N) = \Ext_{\sfMod_{U(L)}}^\bullet(M,N) \quad \text{for all $M,N \in \sfsMod_{U(L)}$.}
	\end{equation}
We will simply write $\Ext_{U(L)}^\bullet(M,N)$ for the common value in \eqref{eq:U(L)Extsmod=Extmod}.

Let $\varphi: k \to k$ be the Frobenius map, $\lambda \mapsto \lambda^2$, and given a $k$-vector space $V$, let $V^{(1)} = k \otimes_{\varphi} V$ denotes its Frobenius twist. Let $\Lone^*[2]$ be the superspace $\Lone^*$ concentrated in $\N$-degree $2$. The squaring operation on $S(\Lone^*)$ induces a homomorphism of graded algebras $S(\Lone^*[2])^{(1)} \to S(\Lone^*)$, $\lambda \otimes_{\varphi} z \mapsto \lambda \cdot z^2$. This is not a homomorphism of graded superalgebras, because the odd superdegree generators in $S(\Lone^*[2])^{(1)}$ are mapped to even superdegree elements in $S(\Lone^*)$. Composing with the inclusion $S(\Lone^*) \hookrightarrow \Lambda(\Lzero^*) \otimes S(\Lone^*) \cong \Cbul(L,k)$, one gets a homomorphism of graded algebras $\wh{\varphi}: S(\Lone^*[2])^{(1)} \to \Cbul(L,k)$. Since $\partial$ acts by derivations on the ring $\Cbul(L,k)$, it follows that the image of $\wh{\varphi}$ consists of cocycles in $\Cbul(L,k)$, and hence $\wh{\varphi}$ induces a graded algebra map
	\begin{equation} \label{eq:varphi}
	\varphi: S(\Lone^*[2])^{(1)} \to \Hbul(L,k),
	\end{equation}
which by abuse of notation we also denote $\varphi$. The map $\varphi$ is natural with respect to $L$, i.e., with respect to even homomorphisms $L \to L'$ of Lie superalgebras.

\begin{theorem} \label{theorem:cfg}
Let $L$ be a finite-dimensional Lie superalgebra over $k$, and let $M$ and $N$ be finite-dimensional $U(L)$-supermodules. Then $\Hbul(L,k)$ is finite over the image of $\varphi$, and $\Ext_L^\bullet(M,N) \cong \Ext_L^\bullet(k,\Hom_k(M,N))$ is finite as either a left or right $\Hbul(L,k)$-module.
\end{theorem}

\begin{proof}
Same as for $p \geq 3$; see \cite{Drupieski:2013c}.
\end{proof}

In fact, a more general version of Theorem \ref{theorem:cfg} is true:

\begin{theorem} \label{theorem:cfg-stronger}
Let $L$ be a finite-dimensional Lie superalgebra over $k$. Let $M$ and $N$ be arbitrary (i.e., not necessarily $\Z_2$-graded) $U(L)$-modules, with $M$ finite-dimensional and $N$ finitely-generated over $U(L)$. Then $\Ext_{U(L)}^\bullet(M,N)$ is finite under the cup product action of $\Hbul(L,k)$. 
\end{theorem}

\begin{proof}
This can be shown in precisely the same manner as in \cite[\S\S2.2--2.3]{Drupieski:2021b}. For the arguments presented there, it is immaterial whether or not the modules are $\Z_2$-graded. The exterior algebra $\Lambda(\wtg_1)$ is a genuine cocommutative Hopf algebra in characteristic $2$, so in the argument at the end of the proof of \cite[Lemma 2.3.1]{Drupieski:2021b}, one can appeal to the results of \cite{Friedlander:1997} rather than \cite[Theorem 3.2.4]{Drupieski:2013c} to argue that the $E_2$-page of \cite[(2.3.5)]{Drupieski:2021b} is a finite module over the cohomology ring $\Hbul(\Lambda(\wtg_1),k)$.
\end{proof}

\subsection{Support varieties}

\emph{For the rest of this section, assume that the field $k$ is algebraically closed.} Taking $A = U(L)$ in \eqref{eq:A-spectrum}, \eqref{eq:A-relative-support}, and \eqref{eq:A-support-variety}, and using the maximal ideal spectrum rather than the prime ideal spectrum, we get the definition of the cohomological spectrum $\abs{U(L)}$, and the definition of (relative) support varieties for $U(L)$-supermodules. Given $M,N \in \sfsMod_{U(L)}$, let $J_{U(L)}(M,N) = \varphi^{-1}(I_{U(L)}(M,N))$, the inverse image of $I_{U(L)}(M,N)$ under the map \eqref{eq:varphi}. Then $J_{U(L)}(M,N)$ is homogeneous with respect to the $\N$-grading on $S(\Lone^*[2])^{(1)}$. Set
	\[
	\Chi_L(M,N) = \Max \left( S(\Lone^*[2])^{(1)} / J_L(M,N) \right).
	\]
and set $\Chi_L(M) = \Chi_L(M,M)$.

\begin{lemma} \label{lemma:relative-homeomorphism}
The graded algebra homomorphism $\varphi: S(\Lone^*[2])^{(1)} \to \Hbul(L,k)$ induces for each pair of finite-dimensional $U(L)$-supermodules $M$ and $N$ a homeomorphism
	\[
	\varphi_{M,N}^*: \abs{U(L)}_{(M,N)} \simeq \Chi_L(M,N).
	\]
\end{lemma}

\begin{proof}
Same as for $p \geq 3$; see \cite[\S4.2]{Drupieski:2019a}.
\end{proof}

If $X$ is an affine variety with coordinate algebra $k[X]$, we write $X^{(1)}$ for the affine variety with coordinate algebra $k[X^{(1)}] = k[X]^{(1)}$. If $Y = X^{(1)}$, then we may write $X = Y^{(-1)}$.

\begin{theorem} \label{theorem:spectrum-odd-nullcone}
Let $L$ be a finite-dimensional Lie superalgebra over an algebraically closed field $k$ of characteristic $p=2$. Then
	\[
	\Chi_L(k)^{(-1)} = \calN_{\odd}(L) := \set{ x \in \Lone : q(x) = 0}.
	\]
\end{theorem}

\begin{proof}
Essentially the same as for $p \geq 3$; see the proof of \cite[Theorem 4.2.4]{Drupieski:2019a}.
\end{proof}

We call $\calN_{\odd}(L)$ the \emph{odd nullcone} of $L$. For $p \neq 2$, one has $\calN_{odd}(L) = \set{ x \in \Lone : [x,x] = 0}$.

\begin{theorem} \label{theorem:support=rank}
Let $L$ be a finite-dimensional Lie superalgebra over an algebraically closed field $k$ of characteristic $p=2$, and let $M$ be a finite-dimensional $L$-supermodule. Then
	\begin{equation} \label{eq:support-equals-rank}
	\Chi_L(M)^{(-1)} = \set{x \in \calN_{\odd}(L) : M|_{\subgrp{x}} \text{ is not free}} \cup \set{0}.
	\end{equation}
\end{theorem}

In \eqref{eq:support-equals-rank}, $M|_{\subgrp{x}}$ denotes the restriction of $M$ to the $k$-subalgebra of $U(L)$ generated by $x$. If $x \neq 0$, this subalgebra has the form $k[x]/(x^2)$, while for $x=0$ the subalgebra is just the field $k$. The set on the right-hand side of \eqref{eq:support-equals-rank} corresponds to the \emph{associated variety} originally considered in characteristic zero by Duflo and Serganova; see \cite{Duflo:2005,Gorelik:2022}.

\begin{proof}[Proof of Theorem 4.6]
Let $\Chi_L'(M)^{(-1)} = \set{ x \in \calN_{\odd}(L) : M|_{\subgrp{x}} \text{ is not free}} \cup \set{0}$ be the `rank variety' defined by the right-hand side of \eqref{eq:support-equals-rank}. For convenience, we will henceforth omit the $(-1)$ superscript. First, using naturality with respect to $L$, one can argue as in the proof of \cite[Proposition 4.3.1]{Drupieski:2019a} to show that $\Chi_L'(M) \subseteq \Chi_L(M)$. To establish the reverse inclusion, we relate support varieties for $U(L)$ to support varieties for $V(L)$ by way of the quotient map $\pi: U(L) \twoheadrightarrow V(L)$. The following argument replaces, for $p=2$, the arguments given in \cite[\S\S3.2--3.3]{Drupieski:2021b} for $p \geq 3$.

Consider the increasing monomial-length filtrations on $U(L)$ and $V(L)$, and the corresponding associated graded algebras
	\begin{gather*}
	\gr U(L) \cong S(\Lzero) \otimes \Lambda(\Lone) \quad \text{and} \\
	\gr V(L) \cong S(\Lzero)/ (x^2 : x \in \Lzero) \otimes \Lambda(\Lone) \cong \Lambda(\Lzero) \otimes \Lambda(\Lone) \cong \Lambda(L).
	\end{gather*}
The algebras $\gr U(L)$ and $\gr V(L)$ are $\N$-graded, with the degree-$0$ component equal to the field $k$, and the degree-$1$ component identifying with the superspace $L$. Now in the same manner described in \cite[\S3.4]{Drupieski:2013c} (cf.\ also \cite[\S2.3]{Drupieski:2021b}), the monomial-length filtrations on $U(L)$ and $V(L)$ induce decreasing filtrations on the corresponding bar complexes for each algebra, and the decreasing filtrations on the bar complexes then give rise to spectral sequences of algebras
	\begin{align}
	{}^VE_1^{i,j} &= \opH^{i+j}(\gr V(L),k)_{-i} \Rightarrow \Hbul(V(L),k), \quad \text{and}  \label{eq:VLspecseq} \\
	{}^UE_1^{i,j} &= \opH^{i+j}(\gr U(L),k)_{-i} \Rightarrow \Hbul(U(L),k). \label{eq:ULspecseq}
	\end{align}
The non-negative gradings on $\gr U(L)$ and $\gr V(L)$ induce non-positive gradings on their cohomology rings, which are indicated by the subscript $-i$ in \eqref{eq:ULspecseq} and \eqref{eq:VLspecseq}.

By the descriptions of $\gr U(L)$ and $\gr V(L)$ given above, one has
	\begin{align*}
	\Hbul( \gr V(L),k) &\cong S(L^*), \quad \text{and} \\
	\Hbul( \gr U(L),k) &\cong \Lambda(\Lzero^*) \otimes S(\Lone^*) \cong S(L^*)/(x^2 : x \in \Lzero^*).
	\end{align*}
In both cases, the superspace $L^* = \Lzero^* \oplus \Lone^*$ is concentrated in the $E_1^{1,0}$ term of the spectral sequence. By multiplicativity, this implies that both spectral sequences are concentrated in the row $j=0$, and consequently the abutment is equal to the cohomology of the $E_1$-page. Since the differential of each spectral sequence acts by algebra derivations, it follows that the square of any element in $S(L^*)$ is a cocycle in ${}^VE_1^{\bullet,0}$. As a complex, ${}^U E_1^{\bullet,0}$ is simply the Koszul cochain complex described above in \eqref{eq:Koszulcochain}. The quotient map $\pi: U(L) \twoheadrightarrow V(L)$ induces a morphism of spectral sequences ${}^VE \to {}^UE$. On the abutment this is just the inflation map induced by $\pi$, while on the $E_1$-page this is the quotient homomorphism $S(L^*) \twoheadrightarrow S(L^*)/(x^2 : x \in \Lzero^*)$.

From the observations of the preceding paragraph, we get the commutative diagram
	\begin{equation} \label{eq:restricted-to-nonrestricted}
	\vcenter{\xymatrix{
	S(L^*[2])^{(1)} \ar@{->}[r]^{\phi} \ar@{->}[d]^{\res} & \Hbul(V(L),k) \ar@{->}[d]^{\pi^*} \\
	S(\Lone^*[2])^{(1)} \ar@{->}[r]^{\varphi} & \Hbul(U(L),k).
	}}
	\end{equation}
Here $\varphi$ is the homomorphism of \eqref{eq:varphi}, and $\phi$ is the analogous map for $\Hbul(V(L),k)$ that arises from the graded $k$-algebra homomorphism $S(L^*[2])^{(1)} \to S(L^*) \cong {}^VE_1^{\bullet,0}$ induced by the squaring operation on $S(L^*)$. Commutativity of \eqref{eq:restricted-to-nonrestricted} then implies commutativity of the diagram
	\begin{equation} \label{eq:pivarietydiagram}
	\vcenter{\xymatrix{
	L^{(1)} \ar@{<-}[rr]^-{\phi^*} \ar@{<-_)}[d]^{\res^*} & & \abs{V(L)} \ar@{<-}[d]^{\pi_*} \ar@{<-^)}[r] & \abs{V(L)}_M \ar@{<-}[d]^{\pi_*} \\
	\Lone^{(1)} \ar@{<-^)}[r] & \Chi_L(k) \ar@{<-}[r]_{\simeq} & \abs{U(L)} \ar@{<-^)}[r] & \abs{U(\wtg)}_M.
	}}
	\end{equation}
By the classical theory (cf.\ \cite{Friedlander:1986b}),
	\[
	\phi^*(\abs{V(L)}_M) = \set{ z \in L : z^{\set{2}} = 0 \text{ and } M|_{\subgrp{z}} \text{ is not free}} \cup \set{0}.
	\]
Then commutativity of the diagram implies that
	\[
	\Chi_L(M) \subseteq \set{z \in \calN_{\odd}(L) : M|_{\subgrp{z}} \text{ is not free}} \cup \set{0} = \Chi_L'(M).
	\]
Then $\Chi_L'(M) \subseteq \Chi_L(M)$ and $\Chi_L(M) \subseteq \Chi_L'(M)$, so $\Chi_L'(M) = \Chi_L(M)$.
\end{proof}

\begin{corollary}[Tensor Product Property] \label{cor:tensor-product-property-U(L)}
Let $L$ be a finite-dimensional Lie superalgebra over an algebraically closed field $k$ of characteristic $p=2$, and let $M$ and $N$ be finite-dimensional $L$-super\-modules. Then
	\[
	\abs{U(L)}_{M \otimes N} = \abs{U(L)}_M \cap \abs{U(L)}_N.
	\]
\end{corollary}

\begin{proof}
Follows immediately from the theorem, in the same manner as Proposition~\ref{prop:tensor-product-property}.
\end{proof}

\begin{theorem}\label{T:DetectFiniteProjDimension}
Let $L$ be a finite-dimensional Lie superalgebra over an algebraically closed field $k$ of characteristic $p=2$. Let $M$ be a finite-dimensional $L$-supermodule. Then
	\[
	\set{ x \in \calN_{\odd}(L) : M|_{\subgrp{x}} \text{ is not free}} = \emptyset \quad \text{if and only if} \quad \projdim_{U(L)}(M) < \infty.
	\]
\end{theorem}

\begin{proof}
The `if' portion of the statement is an immediate consequence of the identification \eqref{eq:support-equals-rank}. For the `only if' portion of statement, the line of reasoning is essentially the same as that for the proof of \cite[Theorem 3.4.2]{Drupieski:2021b}. First, one argues as in the proof of \cite[Lemma 2.4.1]{Drupieski:2021b} that $U(L)$ is free of finite rank over a central (purely even) polynomial subalgebra $\calO$, and hence $U(L)$ is a Noether $\calO$-algebra in the terminology of \cite[Definition A.1.1]{Drupieski:2021b}. Next, one uses Theorems \ref{theorem:support=rank} and \ref{theorem:cfg-stronger} to deduce for each finitely-generated $U(L)$-module $N$ that $\Ext_{U(L)}^i(M,N) = 0$ for all $i \gg 0$. By \cite[Theorem A.1.2]{Drupieski:2021b}, $M$ has finite projective dimension in $\sfMod_{U(L)}$. By \eqref{eq:U(L)Extsmod=Extmod}, this implies that $M$ has finite projective dimension in $\sfsMod_{U(L)}$.
\end{proof}

\begin{remark} \label{remark:Bogvad}
An immediate consequence of Theorem \ref{T:DetectFiniteProjDimension} is that
	\begin{center}
	$\calN_{\odd}(L) = \set{0}$ \quad if and only if \quad $U(L)$ is of finite global dimension.
	\end{center}
On the one hand, if the global dimension of $U(L)$ is finite, then $\Hbul(U(L),k)$ is a finite-dimensional $k$-algebra, and hence $\calN_{\odd}(L) \simeq \abs{U(L)} = \set{0}$. On the other hand, if $\calN_{\odd}(L) = \set{0}$, then $\projdim_{U(L)}(k) < \infty$. Since $\projdim_{U(L)}(M) \leq \projdim_{U(L)}(k)$ for any $U(L)$-supermodule $M$, this implies that $U(L)$ has finite global dimension. This extends to characteristic $2$ a result originally due to B{\o}gvad \cite{Bo-gvad:1984} in characteristic zero (see \cite[Theorem 17.1.2]{Musson:2012} for a general statement in characteristic zero) and due to the authors in odd characteristic \cite[Corollary 3.4.3]{Drupieski:2021b}.
\end{remark}

In any positive characteristic it would be interesting to find a tensor triangulated category whose thick tensor ideals are governed by $\Proj_s (\Hbul(U(L),k))$ in a manner analogous to \eqref{eq:Balmer-spec-homeo-s} and \eqref{eq:BIKP-Theorem10.3-s}. Instead of quotienting out the projective supermodules, as one does for the stable category, one would evidently need to factor out the supermodules of finite projective dimension, since these are the modules on which the support theory for $U(L)$ vanishes. This suggests a version of the singularity category will be the appropriate setting to consider these questions.
	
\makeatletter
\renewcommand*{\@biblabel}[1]{\hfill#1.}
\makeatother

\bibliographystyle{eprintamsplain}
\bibliography{lie-superalgebras-in-characteristic-2}

\end{document}